\documentclass[times]{article}

\usepackage{graphicx}
\usepackage{amssymb}
\usepackage{amsmath}     
\usepackage[usenames, dvipsnames]{color}

\usepackage{fixmath}
\usepackage{multirow}
\usepackage{pgfplots}

\usepackage{tikz}
\usepackage[labelsep=quad,indention=10pt]{caption}
\usepackage{tabu}
\usepackage{float}
\usepackage{subfig}
\usepackage{booktabs}
\usepackage{siunitx}

\usepackage{moreverb}

\usepackage{undertilde}
\usepackage{algorithmic}
\usepackage[lined,boxed]{algorithm2e}

\newcommand{\calV}{\mathcal{V}}
\newtheorem{remark}{Remark}

\newcommand{\bm}[1]{\text{\boldmath $#1$\unboldmath}}

\newcommand{\nS} {\ensuremath{\texttt{n}_{\texttt{S}}}}
\newcommand{\nSL} {\ensuremath{\widetilde{\texttt{n}}_{\texttt{S}}}}
\newcommand{\nQ} {\ensuremath{\texttt{n}_{\texttt{Q}}}}
\newcommand{\nd} {\ensuremath{\texttt{d}}}
\newcommand{\np} {\ensuremath{\texttt{n}_{\texttt{P}}}}

\newcommand\RR{\leavevmode\hbox{$\rm I\!R$}}

\newcommand\bff{\bm{f}}

\newcommand\bn{\bm{n}}
\newcommand\bg{\bm{g}}

\newcommand\bmu{\bm{\mu}}

\newcommand{\TT}{\textsf{T}}
\newcommand\kL{\tilde k}
\newcommand\bK{\bm{K}}
\newcommand\bX{\bm{X}}
\newcommand\bXL{\widetilde{\bX}}
\newcommand\bU{\bm{U}}
\newcommand\bUL{\widetilde{\bU}}
\newcommand\bV{\bm{V}}
\newcommand\bVL{\widetilde{\bV}}
\newcommand\bB{\bm{B}}
\newcommand\bBL{\widetilde{\bB}}
\newcommand\bSigma{\bm{\Sigma}}
\newcommand\bSigmaL{\widetilde{\bSigma}}
\newcommand\bx{\bm{x}}
\newcommand\bu{\bm{u}}
\newcommand\bv{\bm{v}}

\newcommand\bG{\bm{G}}

\newcommand\balpha{\bm{\alpha}}
\newcommand\RB{_{\text{RB}}}
\newcommand\POD{_{\text{POD}}}

\newcommand\bz{\bm{z}}
\newcommand\bw{\bm{w}}
\newcommand\bZ{\bm{Z}}

\newcommand\BibTeX{{\rmfamily B\kern-.05em \textsc{i\kern-.025em b}\kern-.08em
T\kern-.1667em\lower.7ex\hbox{E}\kern-.125emX}}

\begin{document}

\title{Nonlinear dimensionality reduction for parametric problems: a kernel Proper Orthogonal Decomposition (kPOD)}
\author{Pedro D\'iez(1,2), Alba Muix\'i(2), Sergio Zlotnik(1,2),\\ Alberto Garc\'ia-Gonz\'alez(1)\\ \\
$1$- Laboratori de C\`alcul Num\`eric, E.T.S. de Ingenier\'ia de Caminos,\\ Universitat Polit\`ecnica de Catalunya -- BarcelonaTech\\
$2$- International Centre for Numerical\\ Methods in Engineering, CIMNE, Barcelona}

\maketitle

\begin{abstract}
Reduced-order models are essential tools to deal with parametric problems in the context of optimization, uncertainty quantification, or control and inverse problems. The set of parametric solutions lies in a low-dimensional manifold (with dimension equal to the number of independent parameters) embedded in a large-dimensional space (dimension equal to the number of degrees of freedom of the full-order discrete model). A posteriori model reduction is based on constructing a basis from a family of snapshots (solutions of the full-order model computed offline), and then use this new basis to solve the subsequent instances online. Proper Orthogonal Decomposition (POD) reduces the problem into a linear subspace of lower dimension, eliminating redundancies in the family of snapshots. The strategy proposed here is to use a nonlinear dimensionality reduction technique, namely the kernel Principal Component Analysis (kPCA), in order to find a nonlinear manifold, with an expected much lower dimension, and to solve the problem in this low-dimensional manifold. Guided by this paradigm, the methodology devised here introduces different novel ideas, namely: 1) characterizing the nonlinear manifold using local tangent spaces, where the reduced-order problem is linear and based on the neighbouring snapshots, 2) the approximation space is enriched with the cross-products of the snapshots, introducing a quadratic description, 3) the kernel for kPCA is defined ad-hoc, based on physical considerations, and 4) the iterations in the reduced-dimensional space are performed using an algorithm based on a Delaunay tessellation of the cloud of snapshots in the reduced space. The resulting computational strategy is performing outstandingly in the numerical tests, alleviating many of the problems associated with POD and improving the numerical accuracy.

\end{abstract}

\maketitle

\section{Introduction}


Solving parametric problems is an important challenge in computational mechanics. Usually, the set of parametric solutions lies in a manifold of low dimension (equal to the number of independent free parameters, that is the intrinsic dimension of the problem) inside a large-dimensional Euclidean space  (typically, the dimension is the number of degrees of freedom of the discretization). Reduced-Order Models (ROM) are devised to decrease the computational complexity, curtailing the number of degrees of freedom of the full-order problem and enabling real-time solutions of problems requiring multiple queries to the model. The methodology presented in this paper is classified as an \emph{a posteriori} ROM. That is, it is based on sampling the parametric space, computing the corresponding solutions offline, and devising a strategy to drastically diminish the cost of the subsequent online computations (corresponding to different values of the parameters). Alternative \emph{a priori}  approaches construct the ROM offline, explicitly accounting for the parametric dependence with no arbitrary selection of the snapshots \cite{Chinesta-Keunings-Leygue,Chinesta-CLBACGAAH:13,PD-DZGH-18,PD-DZGH-20}: then, the subsequent computations are just functional evaluations.

The most popular a posteriori ROM is Proper Orthogonal Decomposition (POD). POD is extensively used in the framework of dimensionality reduction of boundary value problems to identify linear spaces of lower dimension including the manifold of solutions \cite{Patera-Rozza:07,quarteroni2014reduced,rozza2008reduced}.  The linear space is identified by doing a Principal Component Analysis (PCA) on a set of solutions of the parametric problem (snapshots). Thus, the number of degrees of freedom is readily reduced to the dimension of this linear space using a Reduced Basis technique.
 
Often, the manifold of solutions is \emph{curved}, and therefore the linear subspace embedding the manifold has a dimension much larger than the intrinsic dimension of the parametric family of solutions. The kernel PCA (kPCA) \cite{scholkopf1998nonlinear} technique aims at reducing the large dimensional full-order solutions space (also denoted as \emph{input space} for the role it plays for kPCA) into a nonlinear manifold. 
Following this idea, a kernel POD (kPOD) strategy is introduced here to solve the multidimensional problem in the kPCA nonlinear manifold. Contrary to other nonlinear dimensionality reduction techniques (e.g. Locally Linear Embedding, Isomap... \cite{Isomap2000,LLE2000,Bishop2006}), kPCA is an algorithmically simple approach, proposing an explicit and analytical forward mapping from the full-order space to the reduced-order space.

The driving force of this paper is to take advantage of the kPCA nonlinear dimensionality reduction to solve the parametric problem in a space of a reduced dimension as low as possible, ideally equal to the intrinsic dimension of the problem (number of parameters). Devising a viable a posterior ROM methodology based in this notion requires developing novel strategies that are presented along the paper. The major concepts introduced to complement the kPOD strategy are 
\begin{itemize}
\item using local linear spaces generated by the neighbouring snapshots
\item enriching the approximation space with the cross-products of the snapshots, introducing a quadratic description
\item defining an \emph{ad-hoc} kernel for kPCA, based on the nature of the problem to be solved
\item  introducing a novel iterative algorithm based on a Delaunay tessellation of the cloud of snapshots to solve the nonlinear problem in the reduced space 
\end{itemize}
The resulting strategy is interpretable as a \emph{local} POD approximation, in the sense that for every subsequent query, the approximation is based on a small set of neighbouring snapshots. The idea of selecting a local subset of the snapshots, based on heuristic or automatic clustering techniques is exploited by different authors, see for instance \cite{RAPUN20103046,AMSALLEM2012891,HESS2019379,PAGANI2018530}. Here, the local approach is devised together with the exploration of the nonlinearly reduced space, and naturally associated with the kPCA reverse mapping.

 \begin{figure}[h]
	\centering
	\includegraphics[width = 1\textwidth]{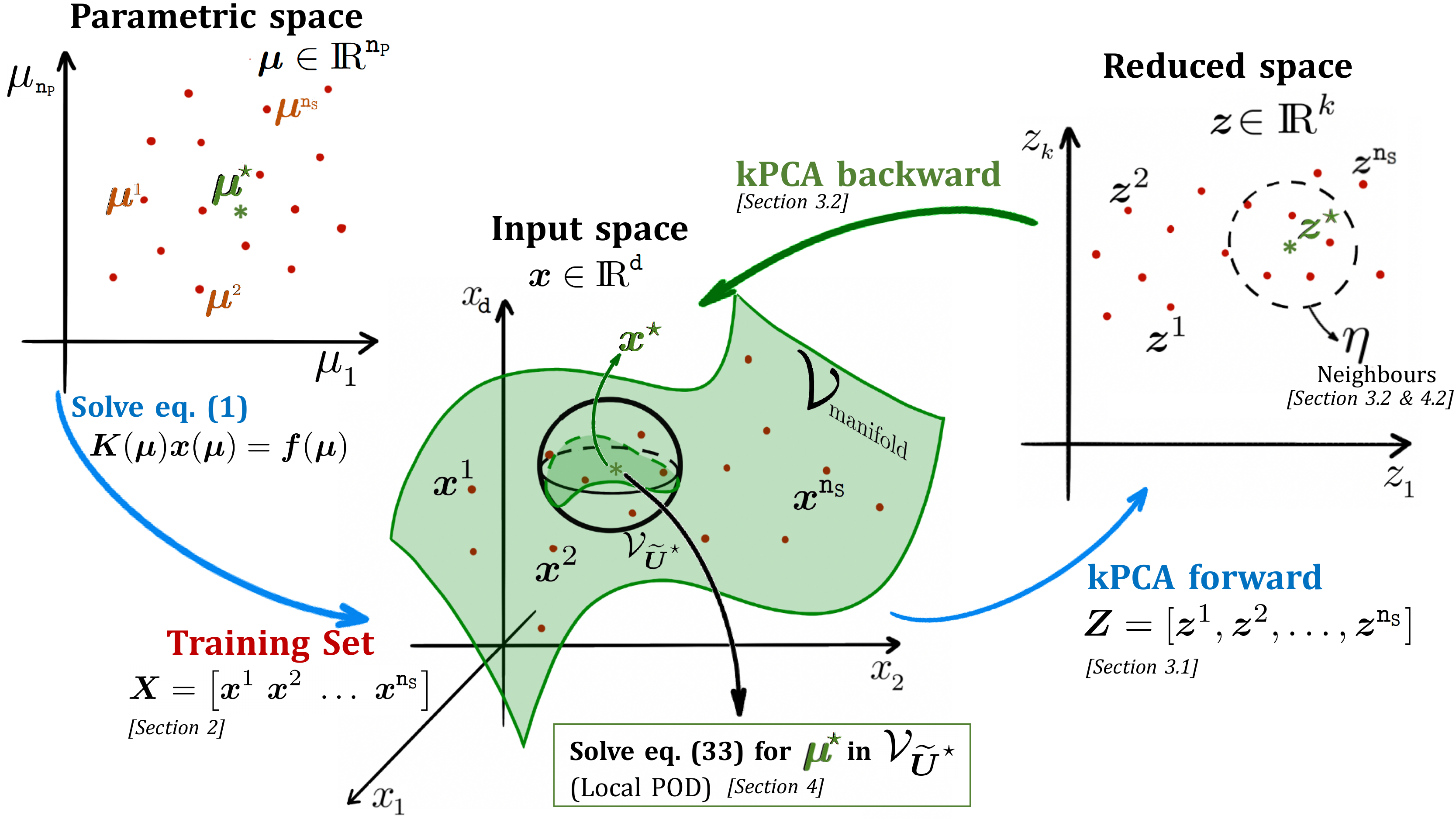}
	\caption{ Illustration of the different steps in the rationale of the proposed methodology.}
	\label{fig:IllustrationPaper}
\end{figure}

The remainder of the paper is structured as follows. First, in section \ref{sec:PbStatement}, the parametric problem is presented, already in discrete form and the standard POD is recalled. Next, in section \ref{sec:kPCA}, the kPCA is briefly described, emphasizing the explicit nature of the forward mapping and the different alternatives for the backward mapping. Also in section \ref{sec:kPCA}, the possibility of enriching the family of snapshots with quadratic terms is introduced.
Section \ref{sec:kPOD} introduces the kPOD concept, with emphasis in the approximation criterion given a local tangent space (subsection \ref{sec:OptiZ}) and the strategy to explore the reduced space  (subsection \ref{sec:exploZ}). Finally, section \ref{sec:NumEx} presents some numerical tests demonstrating the efficiency of the proposed methodology in two 
advection-diffusion problems, and section \ref{sec:Conclu} draws some concluding remarks.

The conceptual structure of the methodology proposed in the paper is summarized in Figure \ref{fig:IllustrationPaper}.  The problem to be solved is a parametric linear system of equations, corresponding to the discretization of some parametric boundary-value or initial-value problem. Solving this problem maps points from the parametric space $\bmu\in\RR^{\np}$ into full-order solutions $\bx(\bmu)\in \RR^{\nd}$ (input space for the dimensionality reduction). Thus, $\nS$ samples $\bmu^{i}$, $i=1,2,\dots,\nS$ (red points in left panel of Figure \ref{fig:IllustrationPaper}) correspond to $\nS$ full-order solutions computed offline in $\RR^{\nd}$ : the snapshots to train the ROM (red dots in the center panel of Figure \ref{fig:IllustrationPaper}). The kPCA technique recalled in section \ref{sec:kPCA} defines the reduced dimension $k\ll \nd$ and maps the training set into the reduced snapshots $\bz^{i}$, $i=1,2,\dots,\nS$ (red dots in the right panel of Figure \ref{fig:IllustrationPaper}, the reduced space). Thus, the kPOD presented in section \ref{sec:kPOD} introduces a methodology to compute the solution associated with a new value $\bmu^{\star}$. This requires exploring the reduced space seeking the optimal $\bz^{\star}$ and the corresponding $\bx^{\star}$ that is best possible approximation to $\bx(\bmu^{\star})$. It is worth noting that each candidate $\bz^{\star}$ defines a local linear space $\calV_{\bUL^{\star}}$ (represented as a sphere in the center panel of Figure \ref{fig:IllustrationPaper}) where the solution is computationally affordable. 

\section{Problem statement, linear dimensionality reduction and linear Reduced Order Model} \label{sec:PbStatement}

\subsection{Problem statement}
The discretization of some parametric Boundary Value Problem results in the following parametric linear system of equations
\begin{equation}\label{eq:LinearSystem}
\bK(\bmu) \bx(\bmu) = \bff(\bmu) ,
\end{equation}
where the unknown $\bx \in \RR^{\nd}$ corresponds typically to the vector of nodal values in a grid, having $\nd$ degrees of freedom. Input data, both matrix $\bK$ and vector $\bff$, depend on the vector of parameters $\bmu \in \RR^{\np}$ (being $\np$ the parametric dimension). The solution $\bx$ depends also on $\bmu$, as it is explicitly marked in the notation.

Thus, the set of solutions $\bx(\bmu)$ lies in a manifold of dimension $\np$ (at the most) inside $\RR^{\nd}$. In this type of problems the number of parameters is much lower than the number of degrees of freedom of the discretization, that is $\np \ll \nd$. Therefore, it is expected that the computational efficiency would drastically improve if the set of solutions in the manifold is described with a substantial reduction of the number of variables, ideally $\np$ if some intrinsic parametrization of the manifold is available. In practice, if the manifold is described by a number of $k$ variables ($\np\le k \ll \nd$), the  size of system of equations \eqref{eq:LinearSystem} is going to be reduced from $\nd$ to $k$.

\subsection{Linear dimensionality reduction: PCA}\label{sec:PCA} 
The first idea is to characterize the manifold in $\RR^{\nd}$ containing all the solutions $\bx(\bmu)$  as a low-dimensional linear space. Obviously, the manifold is not expected to be linear but it may be contained into a linear subspace of a dimension much lower than $\nd$ (in a way, a sort of bounding box).

In order to characterize the manifold, it is populated by a number of samples or \emph{snapshots}. In practice, the parametric space is sampled selecting $\nS$ points in $\RR^{\np}$, say $\bmu^{i}$ for $i=1,2,\dots,\nS$. The corresponding solutions $\bx^{i}=\bx (\bmu^{i})$ are duly computed solving $\nS$ times the full-order system \eqref{eq:LinearSystem}. These solutions constitute a training set that is going to be used to construct the reduced order model.

The snapshots in the training set are collected in a matrix, namely
\begin{equation}\label{eq:Snapshots}
\bX = \left[ \bx^{1}\,\,  \bx^{2}\,\,  \dots \,\, \bx^{\nS} \right] .
\end{equation}

If the sampling is representative of the parametric space, the linear subspace generated by the samples 
\begin{equation*}
\calV_{\bX} = \text{span} \left<  \bx^{1} , \bx^{2}, \dots , \bx^{\nS} \right> \subset \RR^{\nd}
\end{equation*}
is expected to contain all the solutions, also the ones corresponding to the values of $\bmu$ that have not been sampled. Thus,  for a new value of $\bmu$ the original problem \eqref{eq:LinearSystem} is to be solved for $\bx$ in $\calV_{\bX}$ instead of for $\bx$ in $\RR^{\nd}$. As indicated in section \ref{sec:RBPOD} below, the size of the final system to be solved is the dimension of the reduced subspace where the solution is sought, in this case $\calV_{\bX}$.

In order to guarantee that the sampling is representative, the number of samples $\nS$ is taken as large as possible. This is not only increasing the dimension of $\calV_{\bX}$ (and therefore the size of the system to solve); for large $\nS$ the family of snapshots generating $\calV_{\bX}$ may contain linear dependencies (redundancies) resulting in an ill-conditioned reduced system of equations.

In order to eliminate the redundancies and to further reduce the dimension of the linear subspace where the solution is sought, the well known Principal Component Analysis (PCA) technique is used, based in the Singular Value Decomposition (SVD) of $\bX$. The SVD of $\bX \in\RR^{\nd \times \nS}$ produces two unit matrices $\bU\in\RR^{\nd \times \nd}$ and  $\bV\in\RR^{\nS \times \nS}$ and a diagonal matrix $\bSigma\in\RR^{\nd \times \nS}$ such that
\begin{equation}\label{eq:SVD}
\bX = \bU \bSigma \bV^{\TT} ,
\end{equation}
where, for $\nd > \nS$, matrix $\bSigma$ reads
\begin{equation}\label{eq:bSigma}
\bSigma =
\begin{bmatrix}
\sigma_{1} & 0 & 0& \dots & 0 \\
0 & \sigma_{2} & 0& \dots & 0 \\
0 & 0 & \sigma_{3}& \dots & 0 \\
\vdots & \vdots &\vdots & \ddots & \vdots \\
0 & 0 &0& \dots &  \sigma_{\nS}\\
\\
0 & 0 &0& \dots &  0\\
\vdots & \vdots &\vdots & \vdots & \vdots \\
0 & 0 &0& \dots &  0\\
\end{bmatrix} .
\end{equation}
The singular values are sorted in descendent order, $\sigma_{1} \ge \sigma_{2}\ge \dots \sigma_{\nS} \ge 0$.

The columns of $\bU\in\RR^{\nd \times \nd}$, denoted as $\bu^{1}, \bu^{2},\dots,\bu^{\nd}$, are an orthonormal basis of $\RR^{\nd}$. Noting that, \eqref{eq:SVD} is rewritten as 
\begin{equation}\label{eq:SVD2}
\bU  \bSigma = \bX  \bV    \text{ that is }, 
\left[   \sigma_{1} \bu^{1} \, \sigma_{2} \bu^{2} \,\dots \,\sigma_{\nS} \bu^{\nS} \bm{0} \dots \bm{0}   \right] = \bX  \bV 
\end{equation}
it is clear that, if the dimension of $\calV_{\bX}$ is $\nS$ (that is, the snapshots are linearly independent and $\sigma_{\nS} \neq 0$), then the first $\nS$ columns of $\bU$ are an orthonormal  basis of $\calV_{\bX}$.

If the actual dimension of $\calV_{\bX}$ is lower than $\nS$, the lower singular values are zero. In practice, this is equivalent to the case in which they are negligible with respect to the larger ones. In order to keep only the relevant modes, a tolerance $\varepsilon$ is selected and the lower singular values are neglected,  retaining only the $k$ largest values such that
\begin{equation}\label{eq:CollectedVariance}
\sum_{i=1}^{k} \sigma_{i} \ge (1-\varepsilon) \sum_{i=1}^{\nS} \sigma_{i} .
\end{equation}
Thus, taking $\bU^{\star} =\left[\bu^{1} \, \bu^{2} \dots \bu^{k} \right] \in \RR^{\nd \times k}$ as the matrix of the first $k$ columns of $\bU$, the space
\begin{equation}\label{eq:VUstar}
\calV_{\bU^{\star}} = \text{span} \left<  \bu^{1} , \bu^{2}, \dots, \bu^{k}  \right> \subset \RR^{\nd}
\end{equation}
is an approximation to $\calV_{\bX}$ collecting most of the \emph{information}  (up to the tolerance $\varepsilon$).

For instance, if $\varepsilon=0.01$, taking $k$ eigenvalues such that \eqref{eq:CollectedVariance} is satisfied guaranties that $\calV_{\bU^{\star}}$ is a subspace collecting at least 99\% of the \emph{amount of information} (of the variance of the model, to be precise).

\subsection{Linear Reduced Order Model: Reduced Basis and Proper Orthogonal Decomposition}\label{sec:RBPOD}

The Reduced Basis (RB) idea consists in solving the original problem \eqref{eq:LinearSystem} in a subspace of $\RR^{\nd}$ of much smaller dimension but representative of the solutions to be computed.

The crude RB version is to assume that the family of snapshots $\bX$ is a basis, and therefore to solve \eqref{eq:LinearSystem} taking $\bx \approx \bx_{\RB} = \bX \balpha  \in \calV_{\bX}$, where $\balpha \in \RR^{\nS}$ is the vector of unknowns. Thus, system \eqref{eq:LinearSystem} is transformed into
\begin{equation}\label{eq:RBSystemOfEq}
\left[ \bX^{\TT} \bK(\bmu) \bX  \right] \balpha = \bX^{\TT} \bff(\bmu) ,
\end{equation}
which is a system of $\nS$ equations and $\nS$ unknowns.

\begin{remark}[Residual minimization]\label{rem:Enorm}
	In the case $\bK$ is symmetric positive definite, Eq. \eqref{eq:RBSystemOfEq} is equivalent to find $\balpha$ minimizing
	\begin{equation}\label{eq:discrepancyLinear}
	{\mathcal J}(\balpha)=\Vert \bK(\bmu)^{-1} \bff(\bmu)   - \bX \balpha  \Vert_{E}^{2}
	\end{equation}
	Being the $\Vert \cdot  \Vert_{E}$ the discrete form of the standard energy norm defined as 
	\begin{equation}\label{eq:energyNorm}
	\Vert \bx  \Vert_{E}^{2} = \bx^{\TT} \bK(\bmu) \bx.
	\end{equation}
	Note that minimizing the energy norm is consistent with the underlying approximation criterion for the Galerkin formulation.
	For non-symmetric matrices, this equivalence does not hold but the solution of system \eqref{eq:RBSystemOfEq} is all the 
	same an approximation of \eqref{eq:LinearSystem} in $ \calV_{\bX}$. 
\end{remark}
\vspace{2cm}
\begin{remark}[Minimizing residual in Euclidean norm]\label{rem:L2norm} 
	If, instead, the standard Euclidean norm of the residual is selected, that is $\Vert \bx  \Vert^{2} = \bx^{\TT}  \bx$, 
	and the discrepancy form is defined as 
	\begin{equation}\label{eq:discrepancyLinearEuclidean}
	{\mathcal J}(\balpha)=\Vert \bff(\bmu)   - \bK(\bmu) \bX \balpha  \Vert^{2}
	\end{equation}
	the resulting reduced system of equations is 
	\begin{equation}\label{eq:RBSystemOfEqEuclidean}
	\left[\bX ^{\TT}\bK(\bmu)^{\TT} \bK(\bmu) \bX\right]  \balpha = \bX ^{\TT}\bK(\bmu)^{\TT} \bff(\bmu) ,
	\end{equation}
	which is valid also for non-symmetric $\bK$, provided that the symmetric part of $\bK$ is positive definite.
\end{remark}
\vspace{2cm}
Remarks \ref{rem:Enorm} and \ref{rem:L2norm} are pertinent to realize that the reduced system presented in equation \eqref{eq:RBSystemOfEq} is valid in any case, and consistent with the standard Galerkin strategy if $\bK$ is symmetric positive definite. The reduced system \eqref{eq:RBSystemOfEqEuclidean} uses a Least-Squares criterion to project from the full-order space to the reduced-order space, which is not in agreement with the energy projection used from the continuous space to the full-order discrete space. In that sense, the reduced equation  \eqref{eq:RBSystemOfEq} is preferred to \eqref{eq:RBSystemOfEqEuclidean} because it is consistent with the residual minimization for symmetric matrices, and it is adopted also for the case in which $\bK$ is non-symmetric. Both options were implemented in the examples of section \ref{sec:NumEx}, resulting in solutions which are equivalent for all practical purposes.

The family of snapshots becomes a proper basis of $\calV_{\bX}$ if the snapshots are generated guaranteeing that they do not include redundancies. This is the case when a greedy strategy guided with some error assessment is used to generate the basis, see \cite{Florentin2012116,NPM_PD:16,2020-JGR-OZAD}. If the sampling is performed arbitrarily (or even randomly), the family of snapshots typically includes linear redundancies associated with any exhaustive sampling. This results in a rank-deficient matrix $\bX$ (of rank lower than $\nS$) and therefore a singular matrix in system \eqref{eq:RBSystemOfEq} (singular if there are null singular values or extremely ill-conditioned if there are very small with respect to the largest).

The PCA, described in section \ref{sec:PCA}, provides an orthonormal basis of dimension $k \le \nS$ and the associated approximation space $\calV_{\bU^{\star}}$. Combining this with the RB idea leads to the Proper Orthogonal Decomposition concept, that consists in taking $\bx \approx \bx_{\POD} = \bU^{\star} \bz \in \calV_{\bU^{\star}}$, where $\bz \in \RR^{k}$ is the new vector of unknowns. Analogously as for \eqref{eq:RBSystemOfEq}, the POD reduced system reads
\begin{equation}\label{eq:PODSystemOfEq}
\left[  \bU^{\star \TT} \bK(\bmu)  \bU^{\star}  \right] \bz =  \bU^{\star \TT} \bff(\bmu) ,
\end{equation}
which is a system of $k$ equations and $k$ unknowns. 
The POD has eliminated in $ \bU^{\star}$ the possible redundancies in $\bX$. Thus, the reduced system \eqref{eq:PODSystemOfEq} does not spoil the well-conditioned character of matrix $\bK$ (i.e. if $\bK$ is well-conditioned, then also $\bU^{\star \TT} \bK  \bU^{\star}$ is well-conditioned). Actually, the most disadvantageous degradation of the condition number in the reduced system is controlled by tolerance $\varepsilon$. 
Recall that the PCA is a dimensionality reduction technique representing the vectors $\bx\in\RR^{\nd}$ (assumed to be in the set of parametric solutions) in a subspace $\calV_{\bU^{\star}}$ of dimension $k$. 
Thus, a vector $\bx$ is mapped into a $\bz = \bU^{\star \, \TT} \bx \in \RR^{k}$, which is its $k$-dimensional representation. Conversely, $\bz$ is readily mapped back to $\calV_{\bU^{\star}} \subset \RR^{\nd}$ by computing $\bx^{\star} = \bU^{\star} \bz$.  If the original $\bx$ is already in $\calV_{\bU^{\star}}$, $\bx^{\star} = \bx$. In the general case ($\bx \not\in \calV_{\bU^{\star}}$), the pre-image $\bx^{\star}$  is the projection of $\bx$ into $\calV_{\bU^{\star}}$.

It is important noting that in PCA, both the forward mapping $\bx \mapsto \bz =  \bU^{\star \, \TT } \bx $ and the backward mapping $\bz \mapsto \bx^{\star} =  \bU^{\star} \bz$ are explicitly defined.
In the following, the kPCA technique is described as a generalization of PCA to nonlinear multidimensionality reduction and the difficulties associated to properly define the backward mapping are discussed in detail.

\section{Nonlinear dimensionality reduction with kernel-PCA}\label{sec:kPCA}

Kernel Principal Component Analysis (kPCA) is a nonlinear dimensionality reduction technique, based on applying PCA to a transformed training set, being the transformation defined by a kernel form $\kappa(\cdot,\cdot)$ defined in $\RR^{\nd} \times \RR^{\nd}$, see \cite{garcia2020kernel} for details. The kernel $\kappa(\cdot,\cdot)$ is not necessary bilinear but it has to guarantee that the matrix resulting of evaluating the kernel in all the elements of a basis (see below) is symmetric positive definite (this matrix is intended to stand for a Gram matrix). 

Once a kernel $\kappa$ is selected, the Gram matrix $\bG \in \RR^{\nS \times \nS}$ is introduced as having components
\begin{equation}\label{eq:Gmat}
\left[ \bG \right]_{ij}= \kappa(\bx^{i},\bx^{j}) \text{ for } i,j=1,2,\dots,\nS .
\end{equation}

An application $\bg(\cdot)$ from $\RR^{\nd}$ to $\RR^{\nS}$ is introduced associated with kernel $\kappa$:
\begin{equation}\label{eq:gvec}
\bg(\bx) \text{ is such that } \left[ \bg(\bx) \right]_{i}= \kappa(\bx^{i},\bx) \text{ for } i=1,2,\dots,\nS .
\end{equation}
Note that the image by $\bg(\cdot)$ of the $j$-th snapshot, $\bg(\bx^{j})$, coincides with the $j$-th column of $\bG$,  denoted as $\bg^{j}$.

Being $\bG$ symmetric (recall that the kernel $\kappa$ has to guarantee that $\bG$ is symmetric positive definite), its SVD decomposition reads 
\begin{equation}\label{eq:SVDG}
\bG = \bV \bSigma \bV^{\TT} ,
\end{equation}
where in \eqref{eq:SVDG} the $\nS \times \nS$  square diagonal matrix $\bSigma$ has as diagonal entries the squared singular values  $\sigma_{1}^{2} \ge \sigma_{2}^{2}\ge \dots \sigma_{\nS}^{2} \ge 0$.

\subsection{kPCA forward mapping}
The dimensionality reduction is performed setting a tolerance $\varepsilon$ and selecting the reduced dimension $k$ with the criterion given in \eqref{eq:CollectedVariance}. Matrix $\bV^{\star} \in \RR^{\nS \times k}$ is taken as containing the first $k$ columns of $\bV$. Then, the forward mapping on the reduced space $\RR^{k}$ is readily defined for the snapshots: $\bx^{i} \in \RR^{\nd}$ is mapped into  $\bz^{j}=\bV^{\star \, \TT} {\bg}(\bx^{j})$, see \cite{garcia2020kernel} for details. 
The images of the snapshots in the reduced space are collected in matrix $\bZ= [\bz^1,\bz^2, \dots, \bz^{\nS}]$.

Using the definition of $\bg(\cdot)$ in \eqref{eq:gvec}, an exhaustive forward mapping $F$ applied to every $\bx\in\RR^{\nd}$ reads
\begin{equation}\label{eq:kPCAforward}
\begin{aligned}
&F: &&\mathbb{R}^{\nd} &&\longrightarrow &&\mathbb{R}^k \\
& \  &&\bx		  &&\longmapsto     &&\bz = \bV^{\star \, \TT} {\bg}(\bx)
\end{aligned}
\end{equation}

Note that this forward mapping is explicit, similarly as with PCA.

\subsection{kPCA Backward mapping}\label{sec:kPCAbackward}

For the kPCA technique, the question of the backward mapping (or pre-image, or the determination of $F^{-1}$) is not as simple as for the PCA.

Recall that in PCA, the $k$-dimensional linear manifold in $\RR^{\nd}$ is naturally defined as $\calV_{\bU^{\star}}$ in \eqref{eq:VUstar} and therefore, the backward mapping from $\RR^{k}$ to $\calV_{\bU^{\star}}$, maps $\bz$ into $\bx^{\star}= \bU^{\star} \bz$ .

For nonlinear dimensionality reduction techniques, the definition of the pre-image is not as obvious because  it is required to choose the description of the nonlinear  $k$-dimensional manifold in $\RR^{\nd}$. Often, see \cite{mika1998kernel,garcia2020kernel}, the low dimensional manifold is defined locally as a \emph{tangent} space linearly generated by the neighbouring snapshots. \begin{remark}[Tangent space]\label{rem:TangentSpace} 
The concept of \emph{tangent space} refers, here and in the following, to any local linear approximation of the manifold. The dimension of this linear space is typically larger than the dimension of the manifold (here, equal to $k$). In that sense, it differs from the standard geometrical definition of \emph{tangent space}, that has the same dimension of the manifold: here, the \emph{tangent space} is a locally defined linear space that includes the \emph{geometrical tangent space}.
\end{remark}
\vspace{2cm}

In order to define the inverse of the forward mapping $F$ described in \eqref{eq:kPCAforward}, one aims at mapping back any $\bz\in\mathbb{R}^k$  into some $\bx^{\star} \in  \RR^{\nd}$ such that $\bz=F(\bx^{\star})$. 
Function $F$ defined in \eqref{eq:kPCAforward} is at the best surjective (it browses all the elements in $\RR^{k}$) but it cannot be in any case injective (because different elements $\bx \in  \RR^{\nd}$ have the same $\bz$ image). Thus, the inverse is not properly defined unless the target space is a proper low-dimensional manifold $\calV \subset \RR^{\nd}$. The same applies for standard PCA but in this case the linear subspace $\calV_{\bU^{\star}}$ is implicitly selected as target manifold.

Once the low-dimensional manifold $\calV \subset \RR^{\nd}$ is chosen, the backward mapping is taken as
\begin{equation}\label{eq:kPCAbackward}
\bx^{\star} = F^{-1}(\bz) = \arg\min_{\bx\in\calV } \left\Vert F(\bx) - \bz \right\Vert,
\end{equation}
or, in other words, the pre-image is the element $\bx^{\star} \in \calV$ such that $F(\bx^{\star})$ is equal to (or is the best approximation in $\calV$ to) $\bz$.

Thus, the definition of the backward mapping depends on the selection of both the manifold $\calV$ and the approximation criterion, that is the norm $\left\Vert \cdot \right\Vert$.

\subsubsection{Manifold embedded in the linear subspace $\calV_{\bX}$.}\label{ssec:LinearManifold}
The first idea is to select $\calV \subset \calV_{\bX}$ (see section \ref{sec:PCA}),  that is a subset of the space generated by the snapshots. This is equivalent to say that $\bx\in \calV$ is such that for some $\bw = \left[w_{1} \, w_{2} \, \dots w_{\nS} \right]^{\TT} \in \RR^{\nS}$ 
\begin{equation}\label{eq:linearV}
\bx = \sum_{i=1}^{\nS} \bx^i w_i = \bX \bw ,
\end{equation}
namely $\bx\in \calV$ is a weighted average of the snapshots.

The general approach is to select the weights based in the \textbf{minimization of a discrepancy functional} involving the forward mapping, say
\begin{equation}\label{eq:DiscrepancyW}
{\mathcal J}(\bw)=\left\Vert F\left( \bX \bw \right) - \bz \right\Vert^2
\end{equation}
and finding $\bw$ as 
\begin{equation}\label{eq:MinDiscrepancyW}
\bw = \arg\min_{\bm{\omega}\in\RR^{\nS}} {\mathcal J}(\bm{\omega}).
\end{equation}
Having as much as $\nS$ unknown weights would make the minimization problem ill-conditioned.
This is because different combinations of weights (in particular when they correspond to weights associated with snapshots distant of the point of interest) may result in very similar values of $\bx=\bX \bm{\omega}$ that are transformed by $F$ into indistinguishable points $\bz\in\RR^{\nd}$.

 It is therefore necessary to restrict the number of weights to those related with the closest neighbouring snapshots: that is replacing $\bX$ by $\bXL$ and 
$\nS$ by $\nSL$. Even so, solving the minimization problem \eqref{eq:MinDiscrepancyW} is not easy and is typically associated with a number of numerical difficulties (identify a proper first guess, get trapped in local minima...). Often, a more straightforward criterion is often adopted to compute $\bw$ corresponding to $\bz$.

In practice, the weights are selected locally (only the weights corresponding to closer snapshots in $\RR^{\nd}$ are taken as non-zero). Thus, for a given $\bz$ only the closest neighbours are taken into consideration. The set of indices corresponding to these snapshots is denoted by $\eta$ (the dependence of $\bz$ is not explicit in the notation), the number of elements in $\eta$ is denoted by $\nSL$. Correspondingly, the locally reduced matrix of snapshots is denoted by $\bXL$ (the columns of $\bX$ in $\eta$).  The number of unknowns (weights to compute), $\nSL$,  is significantly lower than $\nS$.

This is to say that, for a given $\bz$, the manifold $\calV$ is defined locally as $\calV \subset \calV_{\bXL}$, and there is some $\bw\in\RR^{\nSL}$ such that $\bx=\bXL \bw$. The linear subspace $\calV_{\bXL}$ of dimension $\nSL$ is seen as a \emph{tangent} space to $\calV$ in the vicinity of $F^{-1}(\bx)$, in the sense of Remark \ref{rem:TangentSpace}.

The most classical strategy to compute the weights is using a \textbf{radial basis interpolation}. This consists in computing the distances from $\bz$ to the images of the snapshots in the reduced space, $d_{i} = \left\Vert \bz-\bz^{i}\right\Vert$, for $i=1,2,\dots,\nS$, and define the weights inversely proportional to the squared distances, i.e. $w_{i} \propto \frac{1}{d_{i}^{2}}$. The constant of proportionality is taken such that the sum of the weights is the unity:
\begin{equation}
C= \sum_{j=1}^{\nS} \frac{1}{d_{j}^{2}} \,\,\text{ and } \,\, w_{i} =\frac{1}{C} \frac{1}{d_{i}^{2}}.
\end{equation}
The selection of the squared distances is somehow arbitrary and may be replaced by any monotonic function of the distances (typically with different exponents).
Generally, the criterion to select the set of neighbouring samples $\eta$ is based on the distances $d_{i}$: $i$ is in $\eta$ if $d_{i}$ is lower than some threshold distance. For the radial basis interpolation, this is equivalent to set to zero the weights below some tolerance.

An alternative to compute the weights $\bw$ is proposed here following a classical idea of \textbf{least-squares fitting}.
Note that typically $k\ll \nS$ and even if only $\nSL$ neighbouring snapshots are selected, $k$ is also lower than $\nSL$. Thus, a linear representation of $\bz$ in terms of the snapshots
\begin{equation}\label{eq:LinearZw}
\bz = \bZ \bw
\end{equation}
is an underdetermined system of $k$ equations for $\nS$ (or $\nSL$) unknowns. The typical least-squares-inspired strategy is to use the More-Penrose pseudoinverse of $\bZ\in\RR^{k\times \nS}$, which is a $\nS\times k$ matrix usually denoted by $\bZ^{\dagger}$. Then, the solution of \eqref{eq:LinearZw} is
\begin{equation}\label{eq:LSw}
\bw = \bZ^{\dagger} \bz.
\end{equation}
Note that the solution proposed in \eqref{eq:LSw} is valid for both overdetermined systems of equations ($k>\nS$, which does not hold here; in this case the equality in \eqref{eq:LinearZw} would be an approximation) and underdetermined systems. In the latter case, the pseudoinverse selects among all possible solutions of  \eqref{eq:LinearZw} the one having the lowest norm.

\subsubsection{Manifold including quadratic terms. }\label{ssec:NonLinearManifold}
The main advantage of kPCA with respect to PCA is the possibility of describing a nonlinear manifold $\calV$ of dimension $k\ll \nd$. Actually, the backward mapping strategies proposed in subsection \ref{ssec:LinearManifold} are built as linear tangent approximations, see equation \eqref{eq:linearV}, constituting a globally \emph{curved} manifold. Nevertheless, the local linear approximation may represent a limitation to capture the local curvature of the \emph{actual} manifold $\calV$.

Thus, a new local approximation scheme is introduced accounting for a quadratic representation. The idea is replacing \eqref{eq:linearV} by
\begin{equation}\label{eq:quadV}
\bx = \sum_{i=1}^{\nS} \bx^i w_i +  \sum_{i=1}^{\nS}\sum_{j=i}^{\nS} \bx^i \odot \bx^j w_{ij} =: \bB \bw ,
\end{equation}
where $\odot$ denotes the Hadamard product (component by component). The new vector of weights $\bw \in \RR^{\nQ}$ is of length $\nQ=\nS + \frac{1}{2} \nS (\nS-1)$ (in practice, $\nS$ is to be replaced by $\nSL$ to compute $\nQ$, accounting for the local character of the approximation), and the basis matrix $\bB\in\RR^{\nd \times \nQ}$ is 
\begin{equation}\label{eq:matB}
\bB=\left[\bx^{1}\,\,\bx^{2}\,\,\dots \,\, \bx^{\nS}\,\,\,\,(\bx^{1}\odot\bx^{1})\,\,\,\,(\bx^{1}\odot\bx^{2})\,\,\,\,\dots\,\,\,\, (\bx^{\nS}\odot\bx^{\nS}) \right]
\end{equation}
or the local counterpart $\bBL$ including only $\bx^{i}$ and the different cross-products $\bx^{i}\odot\bx^{j}$, for $i,j\in \eta$.

The same strategies defined in subsection \ref{ssec:LinearManifold} (viz. minimization of a discrepancy functional and least-squares fitting) are to be used here just replacing matrix $\bX$ (or its  restriction $\bXL$ to the selected neighbouring snapshots) by matrix $\bB$ (discarding also the columns involving to the non selected snapshots).
Note that this is equivalent to take a larger tangent space to $\calV$, including also the quadratic terms. In the case of the least-squares fitting, the matrix $\bZ$ has to be augmented including the images by $F$ of the quadratic terms.

\section{A \MakeLowercase{k}PCA-based ROM: the \MakeLowercase{k}POD}\label{sec:kPOD}
The POD strategy leads to a reduced system of equations \eqref{eq:PODSystemOfEq} by seeking, for a given value of $\bmu^{\star}$,  an approximate solution $\bx(\bmu^{\star}) \approx  \bU^{\star} \bz$ of equation \eqref{eq:LinearSystem} in the linear subspace $\calV_{\bU^{\star}}$.


The kPCA proposes a nonlinear manifold, of reduced dimension $k$ which is expected to be lower than the one proposed by PCA and POD. The idea of devising a kPCA-based ROM is to seek the solution $\bx$ of \eqref{eq:LinearSystem} in the nonlinear manifold $\calV$ containing the pre-images $F^{-1}(\bz)$ of every $\bz\in\RR^{k}$.
The approach we present here,  labelled as kernel Proper Orthogonal Decomposition (kPOD), aims at finding in the reduced space the best approximation of the solution of \eqref{eq:LinearSystem} corresponding to a given $\bmu^{\star}\in\RR^{\np}$.
That is exploring all possible $\bz$ (ranging a space of dimension $k$) to pick up $\bz^{\star}\in\RR^{k}$ such that, 
\begin{equation}\label{eq:kPODSystemOfEq}
\bK(\bmu^{\star}) F^{-1}(\bz^{\star}) \approx \bff(\bmu^{\star}).
\end{equation}
Note that the dependence of $\bz^{\star}$ on $\bmu$ is omitted in the notation for the sake of clarity, being the different alternatives to define the backward mapping $F^{-1}$ discussed in section \ref{sec:kPCAbackward}.

The overdetermined system of equations \eqref{eq:kPODSystemOfEq} (with $\nd$ equations and $k\ll \nd$ unknowns)  is to be solved using some fitting criteria that is expressed in a general fashion as the minimization of a discrepancy functional 
\begin{equation}\label{eq:DiscrepancykPCAz}
{\mathcal J}(\bz)=\left\Vert \bK(\bmu^{\star}) F^{-1}(\bz) - \bff(\bmu^{\star})\right\Vert^{2}
\text{ and }
\bz^{\star} = \arg\min_{\bz\in\RR^{k} } {\mathcal J}(\bz),
\end{equation}
Alternatively, the criterion is stated directly for $\bx^{\star}=F^{-1}(\bz^{\star}) \in \calV$,  defining the discrepancy functional analogous to \eqref{eq:discrepancyLinear}, namely  
\begin{equation}\label{eq:DiscrepancykPCAx}
{\mathcal J}(\bx)=\left\Vert \bK(\bmu^{\star})^{-1} \bff(\bmu^{\star})   - \bx  \right\Vert_{E}^{2}
\text{ and }
\bx^{\star} =  \arg\min_{\bx\in\calV } {\mathcal J}(\bx) .
\end{equation}

According to expression \eqref{eq:DiscrepancykPCAz} the space to explore in the minimization is of dimension $k$, associated with the kPCA reduction, which is expected to be lower than the one provided by POD (the kPCA nonlinear dimensionality reduction is likely closer to the intrinsic dimension of the problem). The price to pay for reducing the $k$ is, in this case, the fact that problem \eqref{eq:DiscrepancykPCAz} is nonlinear while the POD reduced system \eqref{eq:PODSystemOfEq} keeps the linear character of the original problem \eqref{eq:LinearSystem}. The advantage of using formulation \eqref{eq:DiscrepancykPCAx} is that once the set of neighbouring snapshots $\eta$ (see section \ref{ssec:LinearManifold}) is defined, the local neighbourhood of manifold $\calV$ is embedded in a linear space similar to $\calV_{\bXL}$ (or $\calV_{\bBL}$ if the quadratic terms are also included) and the problem to be solved is analogous to the reduced system of the POD, see equation \eqref{eq:PODSystemOfEq}.

\subsection{Optimal solution for a given point $\bz$ in the reduced space}\label{sec:OptiZ}
The given point $\bz$ in the reduced space is associated with a set of neighbouring snapshots, whose indices are collected in $\eta$. Thus, the characterization of the backward mapping $F^{-1}(\cdot)$ is, in practice, a matter of selecting a proper distribution of weights such that $F^{-1}(\bz^{\star})=\bx^{\star}=\bBL\bw$, as described in section \ref{sec:kPCAbackward}. This is equivalent to determining an element in the tangent space $\calV_{\bBL}$. 

In this context, to improve the local approximation properties, it is relevant to center the selected set of samples and generate the quadratic terms correspondingly. Namely, to compute the local average
\begin{equation}\label{eq:localAver1}
\bar\bx = \frac{1}{\nSL} \sum_{i\in \eta} \bx^{i}
\end{equation}
and to express $\bx$ in the tangent space as (including quadratic terms)
\begin{equation}\label{eq:localAver2}
\bx = \bar\bx + \sum_{i\in \eta} (\bx^{i} - \bar\bx) w_{i} 
+  \sum_{i\in \eta}  \sum_{\tiny{\begin{matrix} j \in \eta \\  j\ge i \end{matrix}} }
(\bx^{i} - \bar\bx)\odot  (\bx^{j} - \bar\bx) w_{ij} = \bar\bx + \bBL \bw .
\end{equation}

In order to eliminate possible redundancies in $\bBL\in\RR^{\nd\times \nQ}$, following equation \eqref{eq:SVD}, the SVD is performed and reads $\bBL=\bUL \bSigmaL \bVL^{\TT}$. The matrix is reduced selecting $\kL$ with a criterion similar to \eqref{eq:CollectedVariance}. Then, matrix $\bUL^{\star}\in\RR^{\nd\times \kL}$ containing the first $\kL$ columns of $\bUL$ is used to define the reduced approximation 
\begin{equation}\label{eq:redLocalApprox}
\bx^{\star}= \bar\bx + \bUL^{\star} \bw^{\star} .
\end{equation}

The manifold $\calV$ (locally approximated by the affine subspaces $\{ \bar\bx \} + \calV_{\bBL}$ or $\{ \bar\bx \} +\calV_{\bUL^{\star}}$) contains the solutions of  \eqref{eq:LinearSystem}. The optimal value of $\bw^{\star}$ is therefore readily computed solving the $\kL\times \kL$ linear system analogous to \eqref{eq:PODSystemOfEq} 
\begin{equation}\label{eq:kPODlocalSystem}
\left[  \bUL^{\star \TT} \bK(\bmu^{\star})  \bUL^{\star}  \right] \bw^{\star} =  \bUL^{\star \TT} \left[  \bff(\bmu^{\star}) - \bK(\bmu^{\star}) \bar\bx \right].
\end{equation}

\subsection{Exploration of the reduced space}\label{sec:exploZ}

The strategy proposed in section \ref{sec:OptiZ} procures the solution according to the criterion expressed in \eqref{eq:DiscrepancykPCAx}. The only difference lies in the set $\eta$, associated with the space where $\bx^{\star}$ is selected for the minimization problem:  instead of the nonlinear manifold $\calV$, the minimization is performed in the tangent affine space $\{ \bar\bx \} +\calV_{\bUL^{\star}}$, which depends on the point $\bz$ selected in the reduced space $\RR^{k}$.

As it is clear from \eqref{eq:DiscrepancykPCAz}, the reduced space has to be explored to find the point $\bz$ producing the best solution. The quality of the different approximations is to be assessed using some residual-based discrepancy functional, see \eqref{eq:DiscrepancykPCAz}.

It is also clear from section  \ref{sec:OptiZ} that the role of $\bz$ in the computation of $\bx^{\star}$ is only to define the closest neighbours, that is, the set $\eta$ containing the indices of the $\nSL$ neighbouring snapshots. Accordingly, it makes sense to explore $\RR^{k}$ discretely, browsing the zones in which the neighbouring points in the training set are the same, that is \emph{patches} around each snapshot.

Thus, the strategy proposed here to explore the different $\bz$ boils down to define patches in $\RR^{k}$, associated with the cloud of points in the training set mapped into the reduced space.  These  patches are in fact characterizing the neighbours connected to each snapshot. 

\begin{figure}
	\centering
	\includegraphics[width = 0.99\textwidth]{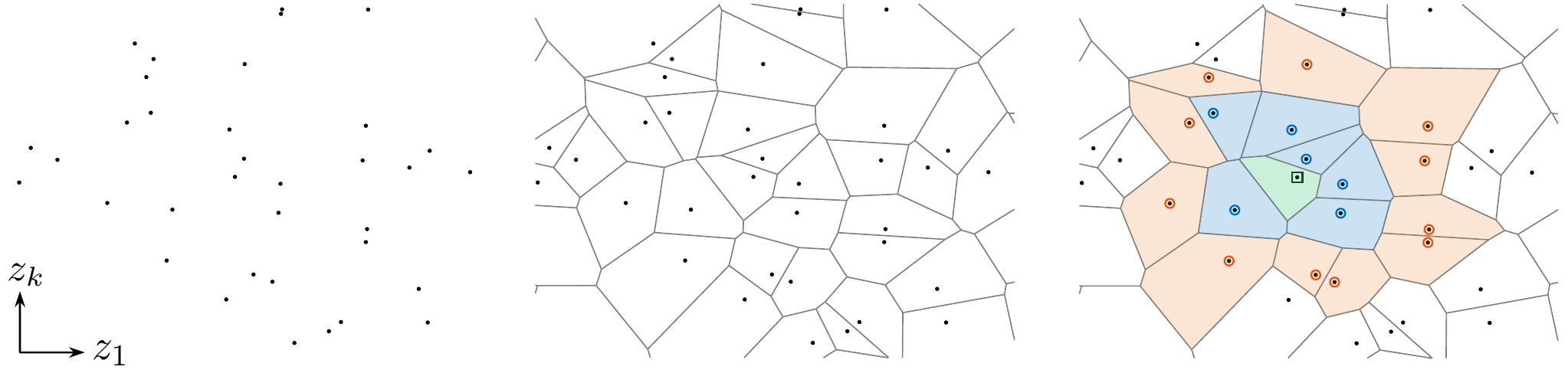}
	\caption{The cloud of snapshots (left) defines a Voronoi tessellation of cells associated with each snapshot (center). The right panel illustrates the levels of connectivity associated to some snapshot $\bz^{i}$ (marked as $\square\!\!\!\!\,\bm{\cdot}\,\,\,$). Connectivity is based on the Voronoi diagram. The first-level patch includes $\bz^{i}$ itself and all the snapshots  whose Voronoi cells are adjacent to it (marked with 
${\color{blue}\bigcirc}\!\!\!\!\!\,\bm{\cdot}\,\,\,$). The second-level patch includes also the cells associated with snapshots marked with 
${\color{red}\bigcirc}\!\!\!\!\!\,\bm{\cdot}\,\,\,$.
}
\label{fig:VoronoiLevels}
\end{figure}

Here, this is readily performed using a generalized Delaunay tessellation (and the equivalent Voronoi diagram) of the training set in $\RR^{k}$. Note that the Voronoi cells $\Omega_{i}\subset \RR^{k}$ are centered in the points $\bz^{i}$ of the training set (snapshots), for $i=1,2,\dots,\nS$. A straightforward criterion is to consider as neighbours of $\bz\in\Omega_{i}$ all the points  connected to $\bz^{i}$ in the Delaunay tessellation, see Figure \ref{fig:VoronoiLevels}. This can be taken with one level of connectivity (being vertices of the same simplex) or two. Figure \ref{fig:VoronoiLevels} illustrates how the cloud of points in $\RR^{k}$ (in the illustration $k=2$) is tessellated establishing both the Voronoi cells and the connectivity collected in the set of indices $\eta$. For a given snapshot $\bz^{i}$ (point marked with $\square\!\!\!\!\,\bm{\cdot}\,\,\,$ in the figure), the set $\eta$ may account for one level of surrounding cells  (those marked with ${\color{blue}\bigcirc}\!\!\!\!\!\,\bm{\cdot}\,\,\,$) or two (adding also those marked with ${\color{red}\bigcirc}\!\!\!\!\!\,\bm{\cdot}\,\,\,$)

The use of the Delaunay algorithm is an excellent option to select the neighbouring points because it is available in many platforms and it is generalized to arbitrary dimensions (no limit in the value of $k$), see \cite{george1998delaunay}.

\begin{figure}
	\centering
	\includegraphics[width = 0.99\textwidth]{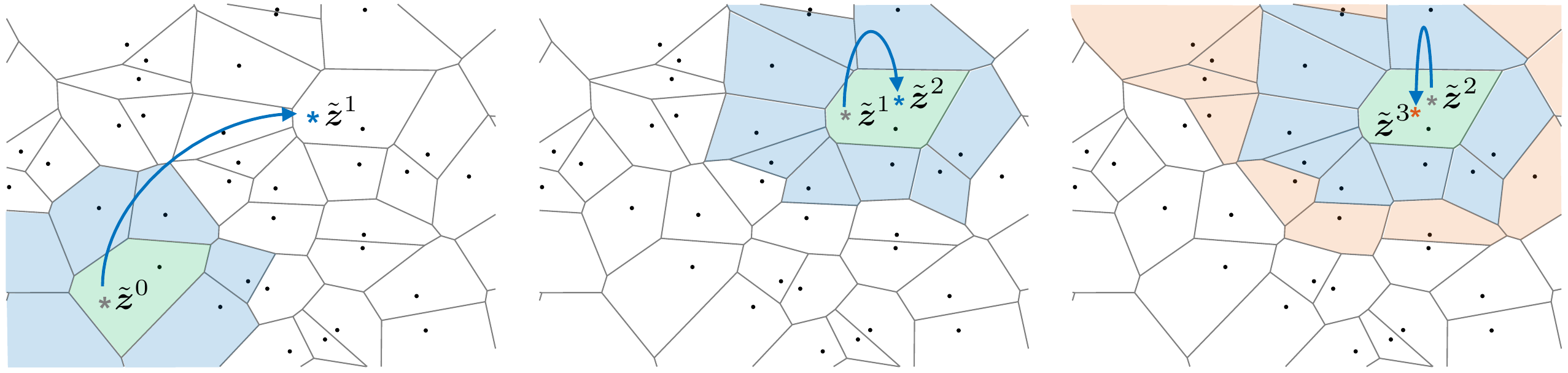}
	\caption{Illustration of the \emph{optimal path} strategy to iteratively explore the $\bz$ space. The patch of the initial guess $\tilde\bz^{0}$ produces a local solution $\tilde\bz^{1}$ using the first connectivity level  (left panel). If $\tilde\bz^{1}$ is not in the same Voronoi cell of $\tilde\bz^{0}$, the next iteration $\bz^{2}$ is computed as the local solution in the patch of $\tilde\bz^{1}$ with one connectivity level  (center panel). If $\tilde\bz^{2}$ is in the same cell $\Omega_{j}$ as $\tilde\bz^{1}$, the next iteration $\tilde\bz^{3}$ is computed increasing the level of connectivity (right panel). The iterative algorithm is stopped if $\tilde\bz^{3}$ is in the same cell as $\tilde\bz^{2}$ (and hence, as $\tilde\bz^{1}$).}
\label{fig:OptimalPathIteration}
\end{figure}

The strategy described above is defining the set $\eta$ associated with any candidate point in the reduced space $\bz$. Once $\eta$ is determined (having chosen if the neighbourhood takes one or two levels of Voronoi cells), the corresponding $\bx^{\star}$ is computed as indicated in section \ref{sec:OptiZ}. The question is whether this solution is acceptable or if there is need of a further iteration in the reduced space. 

The idea presented here, and illustrated in Figure \ref{fig:OptimalPathIteration}, is to devise an iterative algorithm generating a succession of points in $\RR^{k}$, 
$\tilde\bz^{0},\tilde\bz^{1},\dots,\tilde\bz^{\nu},\tilde\bz^{\nu+1},\dots$ The iteration scheme starts with some initial guess $\tilde\bz^{0}$, taking just one level of Voronoi cells around the cell $\Omega_{i}$ where $\tilde\bz^{0}$ is located. This characterizes $\eta$ and therefore allows computing $\bx^{\star}$ associated with $\tilde\bz^{0}$. Then, $\bx^{\star}$ is mapped into the reduced space using the kPCA forward mapping, see equation \eqref{eq:kPCAforward} and  $\tilde\bz^{1}$ is obtained. In the case $\tilde\bz^{1}$ is in a Voronoi cell $\Omega_{j}$ different than the previous ($j\neq i$), the computation is performed again selecting $\eta$ using a one-level Voronoi connectivity around $\Omega_{j}$. The operation is repeated until two successive iterations $\tilde\bz^{\nu}$ and $\tilde\bz^{\nu+1}$ belong to the same Voronoi cell, say $\Omega_{j}$. In order to avoid a possible local stagnation and to eventually increase the accuracy of the proposed solution, in this case the computation is performed again selecting a two-level connectivity around $\Omega_{j}$. If, as it is observed in all the examples in section \ref{sec:NumEx}, the next iteration, say $\tilde\bz^{\nu+2}$ is also in the same cell $\Omega_{j}$, the iterative process is stopped. Else, if $\tilde\bz^{\nu+2}$ is in $\Omega_{\ell}$ with $\ell\neq j$, the iterative process continues. Besides being illustrated in Figure \ref{fig:OptimalPathIteration}, this process is algorithmically described in the next section.

One important point is properly selecting the initial guess $\tilde\bz^{0}$. A very effective possibility is to compute a POD solution and to map it to the reduced space using the standard forward kPCA mapping  \eqref{eq:kPCAforward}.

\subsection{Algorithmic description}\label{sec:Algo}
The proposed methodology is summarized here using an algorithmic description. The process is split into an offline phase and an online phase. The latter corresponds to the iterative procedure described in section \ref{sec:exploZ}. The offline phase includes the computation of the snapshots (full-order problems corresponding to selected points in the parametric space), the kPCA analysis and the pre-process of the training set in the reduced space (Delaunay tessellation). These two phases are synthetically described next and summarized in algorithms \ref{alg:OfflinePhase} and \ref{alg:OnlinePhase}.

The first offline phase consists in 
\begin{enumerate}
	\item sampling the parametric space in points $\bmu^{i}$ for $i=1,2,\dots,\nS$ and compute the corresponding snapshots $\bx^{i}$ (training set $\bX\in\RR^{\nd\times  \nS}$) solving equation \eqref{eq:LinearSystem} $\nS$ times.
	\item performing the kPCA dimensionality reduction: compute matrix $\bG$, perform SVD as in equation \eqref{eq:SVDG}, and find the reduced dimension $k$ and matrix $\bV^{\star}$ based in the tolerance $\varepsilon$.
	\item map the snapshots into $\RR^{k}$ and compute $\bZ\in  \RR^{k \times \nS}$
	\item Build a Delaunay tessellation in $\RR^{k}$ with the reduced snapshots $\bZ$, characterize cells $\Omega_{i}$ and the snapshots connected to them.
\end{enumerate}

\begin{algorithm}[h!] 
	\begin{algorithmic}
		\STATE{\textbf{Input Data:} Sampled parameters $\bmu^{i}$ for $i=1,2,\dots,\nS$\\ \hskip2.2cm 
			matrix $\bK(\mu)$ and vector $\bff(\bmu)$, tolerance $\varepsilon$, kernel $\kappa$.}
		\\ $\,$ \\
		\STATE{Compute $\bx^{i}$ solving $\bK(\bmu^{i}) \, \bx^{i}= \bff(\mu^{i})$  for $i=1,2,\dots,\nS$ }
		\STATE{Collect all $\bx^{i}$ in matrix $\bX$}
		\STATE{Compute matrix $\bG$ using $\kappa$}
		\STATE{Perform SVD: $ \bG = \bV \bSigma \bV^{\TT}$, compute $k$ and $\bV^{\star}$ (using $\varepsilon$)}
		\STATE{Compute matrix of reduced snapshots $\bZ=\bV^{\star \, \TT} \bG$ }
		\STATE{Build a generalized Delaunay Tessellation of $\bZ$ in $\RR^{k}$: $\Omega_{i}$, $i=1,\dots,\nS$}
		\\ $\,$ \\
		\STATE{{\textbf{Output Results:} snapshots $\bX$ and $\bZ$, matrix $\bV^{\star}$, tessellation  $\Omega_{i}$ }}
	\end{algorithmic}
	\caption{Algorithmic description of the offline phase.}\label{alg:OfflinePhase}
\end{algorithm}

Then, for a new value of $\bmu$, the reduced order solution is computed in an online phase by
\begin{enumerate}
	\item loop in $\tilde\bz^{\nu}$, for $\nu=0,1,2,\dots$ until convergence
	\item find the cell $\Omega_{i}$ to where $\tilde\bz^{\nu}$ belongs (centered in $\bz^{i}$), for some $i\in\{1,2,\dots,\nS\}$
	\item select the neighbours based in the Delaunay connectivity and compute local average $\bar\bx$ and the local basis $\bBL$ (including quadratic terms if needed). In order to guarantee the quality of the local approximation, for snapshots located close to the boundary of the $z$-domain where the number of natural neighbours is small, a minimum number of surrounding  snapshots is enforced according the dimension of the reduced space, $k$, and the selected number of levels. 
	\item perform SVD of $\bBL$, find reduced dimension $\kL$, matrix $\bUL^{\star}$ and compute $\bx^{\star}$ solving equation \eqref{eq:kPODlocalSystem} ($\kL \times \kL$ system) and \eqref{eq:redLocalApprox}
	\item compute next iteration with the forward kPCA image of  $\bx^{\star}$,
	$\tilde\bz^{\nu+1}=\bV^{\star \, \TT} {\bg}(\bx^{\star})$
	\item  find the cell $\Omega_{j}$ to where $\tilde\bz^{\nu+1}$ belongs
	\begin{itemize}
\item[ ] If $\Omega_{j}\neq \Omega_{i}$ (and $\Omega_{j}$ has not been visited yet): go to step 3
\item[ ] Else (if  $\Omega_{j} = \Omega_{i}$ or $\Omega_{j}$ has been already visited): perform an extra computation in the $\Omega_{j}$ patch to check stagnation and to increase accuracy
\end{itemize}
	\item keep $\bx^{\star}$ as solution
\end{enumerate}

\begin{algorithm}[h!] 
	\begin{algorithmic}
		\STATE{\textbf{Input Data:} New value of $\bmu$,  matrix $\bK(\bmu)$ and vector $\bff(\bmu)$ \\ \hskip2.15cm 
			tolerance $\varepsilon$, kernel $\kappa$,  snapshots $\bX$ and $\bZ$,  \\ \hskip2.15cm 
			matrix $\bV^{\star}$, tessellation $\Omega_{i}$, $i=1,2,\dots,\nS$}
		\\ $\,$ \\
		\STATE{Initialize: select $i$ in $\{1,2,\dots,\nS \}$}
		\WHILE{Stopping criterion is not satisfied} 
		\STATE{Find the set of neighbouring indices, $\eta\subset \{1,2,\dots,\nS \}$} 
		\STATE{Build $\bXL$ as the columns of $\bX$ in $\eta$}
		\STATE{Compute local average $\bar\bx$ and local basis $\bBL$}
		\STATE{Perform SVD: $ \bBL = \bUL \bSigma \bVL^{\TT}$, compute $\kL$ and $\bUL^{\star}$ (using $\varepsilon$) }
		\STATE{Solve $\left[  \bUL^{\star \TT} \bK(\bmu)  \bUL^{\star}  \right] \bw^{\star} 
		=  \bUL^{\star \TT} \left(\bff(\bmu) - \bK(\bmu) \bar\bx \right) $} 
		\STATE{Compute $\bx^{\star}= \bar\bx + \bUL^{\star} \bw^{\star}$}
		\STATE{Compute $\bz^{\star}=\bV^{\star \, \TT} \bg(\bx^{\star})$}
		\IF{$\bz^{\star}$ is in  $\Omega_{i}$}
		\STATE{Re-compute $\bz^{\star}$ with an additional level of connectivity}
		\IF{ $\bz^{\star}$ continues in $\Omega_{i}$}
		\STATE{Stopping criterion is  satisfied}
		\ELSE
		\STATE{Set new value of $i$ such as $\bz^{\star}\in\Omega_{i}$}
		\ENDIF	
		\ELSE
		\STATE{Set new value of $i$ such as $\bz^{\star}\in\Omega_{i}$}
		\ENDIF
		\ENDWHILE
		\\ $\,$ \\
		\STATE{\textbf{Output Result:} solution $\bx^{\star}$ }
	\end{algorithmic}
	\caption{Algorithmic description of the online phase.}\label{alg:OnlinePhase}
\end{algorithm} 

\section{Numerical examples}\label{sec:NumEx}

This section presents two numerical examples involving a advection-diffusion problem. 
The first example is a transient one-dimensional problem, while the second one solves a steady state problem in a two-dimensional domain. 
In the two examples, the physical properties of the problem are exploited to build an ad hoc kernel function for kPOD, which successfully identifies the intrinsic dimension of the manifold of solutions.

Throughout the section, kPOD includes quadratic terms in the approximation space. 
Errors are measured with the Euclidean norm with respect to a reference solution, which is taken as the full-order FE solution. 

\subsection{Transient advection-diffusion problem in 1D}

\begin{figure}[]
	\centering
	\includegraphics[width=0.65\textwidth]{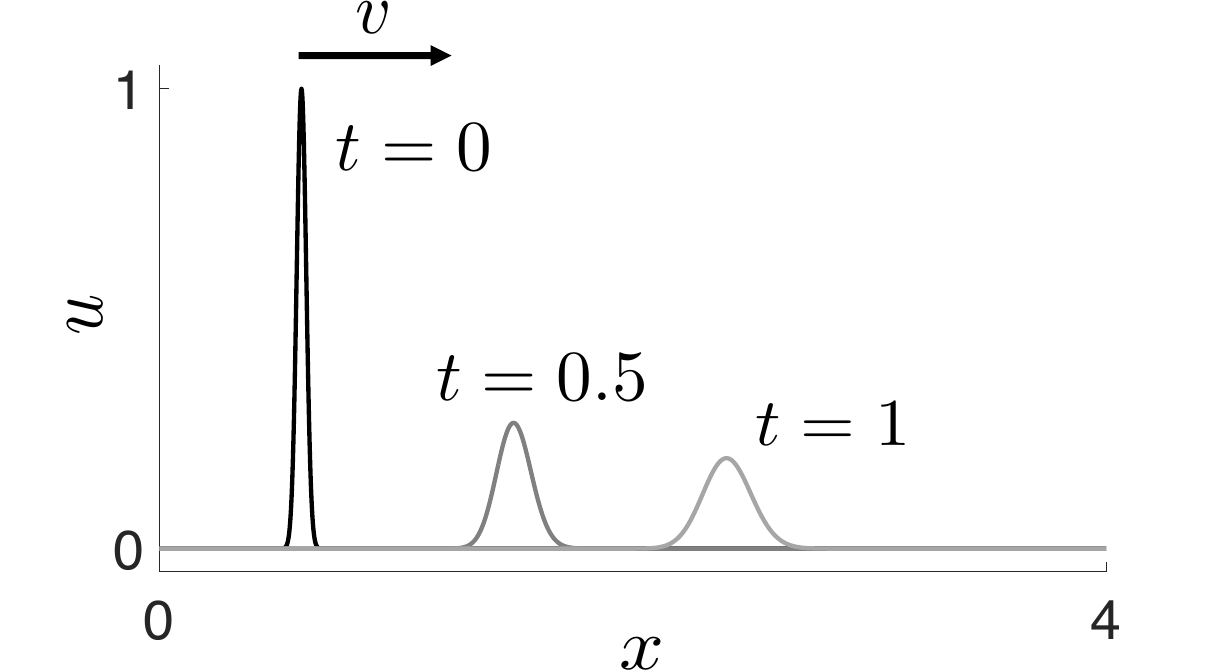}
	\caption{Illustration of the 1D-transient advection-diffusion solution by three representative time snapshots.
	The initial profile $u(x,0)$ (in black) diffuses and propagates to the right as time evolves.}
	\label{fig:test1D-expected}
\end{figure}

We consider the transient boundary value problem
\begin{equation}\label{problem1D}
\left\{
\begin{aligned}
&u_t + v u_x + \nu u_{xx} = 0  &\text{ for } x\in(0,4),\, t\in(0,1.25],\\ 
&u(0,t) = u(4,t) = 0, & \\
& u(x,0) = \exp\left(-\frac{1}{2} \left( \frac{x-0.6}{0.02} \right)^2 \right),&
\end{aligned}
\right.
\end{equation}
where the diffusivity is $\nu = 5\cdot 10^{-3}$ and the advection velocity $v$ takes values in the range $[1,2]$.
For some $v$,  the expected evolution of  a representative solution in time is sketched in Figure \ref{fig:test1D-expected}.

The problem is to be solved when varying two parameters: the time and the advection velocity. 
Nevertheless, due to the small diffusivity, $vt$ is almost constant as the solution evolves in time, and the solutions lie in (or are very close to) a one-dimensional manifold.

The spatial discretization  is uniform with $2000$ intervals, leading to nodal solutions in $\RR^{2001}$. In time, we take $250$ steps with an increment of $\Delta t = 0.005$. Discretizing equation \eqref{problem1D} using centered finite differences in space and Crank-Nicolson in time, the nodal solution $\bu^{n+1}$ at time step $n+1$ is computed by solving a system of the form
\begin{equation}\label{transient-full}
\bm{A}(v) \bu^{n+1} = \bm{D}(v) \bu^{n},
\end{equation}
where $\bu^n$ is the nodal solution from time step $n$.
Accounting for the Dirichlet values in the extremes of the interval, $\bm{A}$ and $\bm{D}$ are matrices of size $1999\times1999$. Our goal is to reduce the size of these systems by means of POD and kPOD, in order to obtain a more efficient computation at every time iteration.

The training set for this problem consists of $501$ snapshots, including the inital condition and the full-order solutions of the problem for $10$ equiespaced velocities in  $[1,2]$, for 50 time steps corresponding to $n = 5, 10, 15 \dots 250$. 

Within the POD formulation, system \eqref{transient-full} becomes
\begin{equation}
\begin{aligned}
&\left[ \bU^{\star \TT} \bm{A}(v) \bU^\star \right] \bm{v}^{n+1} = \bU^{\star \TT} \left[ \bm{D}(v) \bu^{n} - \bm{A}(v) \bar{\bu} \right], \\
&\bu^{n+1} = \bar{\bu} + \bm{v}^{n+1},
\end{aligned}
\end{equation}
where $\bar{\bu}$ is the mean of the snapshots and $\bU^\star$ is used to define the reduced-order approximation as indicated in \eqref{eq:redLocalApprox}. Taking a tolerance $\varepsilon = 10^{-8}$ in the dimensionality reduction criterion \eqref{eq:CollectedVariance}, the resulting system of equations is of size $74\times 74$.

With kPOD, it is necessary to solve, for each time step, as many systems as determined by the optimal path algorithm, described in Section \ref{sec:Algo}. Thus the system for each iteration in  algorithm \ref{alg:OnlinePhase} reads
\begin{equation}
\begin{aligned}
&\left[ \tilde{\bU}^{\star \TT} \bm{A}(v) \tilde{\bU}^\star \right] \bm{w}^{n+1} = \tilde{\bU}^{\star \TT} \left[ \bm{D}(v) \bu^{n} - \bm{A}(v) \bar{\bu} \right], \\
&\bu^{n+1} = \bar{\bu} + \tilde{\bU}^\star \bm{w}^{n+1},
\end{aligned}
\end{equation}
with a $\tilde{k}\times\tilde{k}$ system matrix, $\tilde{\bU}^\star$ and $\tilde{k}$ chosen as described in Section \ref{sec:OptiZ}. The dimension $\tilde{k}$ depends on the particular neighborhood of snaphots under study, also taken with a tolerance $\varepsilon = 10^{-8}$.
Here, the algorithm at time step $n+1$ is initialized at $\tilde{\bz}^0 = F(\bu^n)$. 

The choice of an appropriate kernel function is relevant to guarantee the efficiency and the accuracy of the method. 
Often, the kernel is selected agnostic to the problem, based on the Euclidian distance of the snapshots in the full-order space $\RR^{\nd}
$. Here, the kernel is selected such that the \emph{distance} between samples corresponds to physical considerations, taking advantage of the a priori knowledge on the nature of the solutions. The expected evolution of the solutions is depicted in Figure \ref{fig:test1D-expected}: the initial profile propagates to the right with time, due to advection, and decreases in heigh and increases in width, keeping the area constant, due to diffusion. Thus, the solution $u(x;t,v)$ is mainly characterized by the position of its centroid. In other words, the horizontal and vertical coordinates of the centroid are describing the progress of the initial configuration, accounting for both advection (horizontal position) and diffusion (vertical position). For any $u: [a,b] \rightarrow \mathbb{R}$, the centroid is defined as 
\begin{equation}\label{centroid_definition}
\bm{C}_{[a,b]} (u) = \frac{1}{\int_a^b u(x) dx} \left( \int_a^b x u(x) dx , \int_a^b \frac{u^2(x)}{2} dx\right).
\end{equation}
Thus, a sensible choice for the kernel in this problem results from replacing in a classical Gaussian kernel the Euclidean distance of the solutions by the \emph{centroid distance}, namely
\begin{equation}
\kappa(\bu,\bv) = \exp \left(-\beta \left\| \bm{C}_{[0,4]}(\bu) - \bm{C}_{[0,4]}(\bv) \right\|^2 \right),
\end{equation}
$\bu,\bv$ being two solutions (snapshots) and $\beta$ a constant, here is set to $10^{-4}$. The value of $\beta$ is chosen after numerical experimentation, trying to accumulate a significant amount of information for $k = 1$. In fact, with the proposed kernel the first principal component collects $99.73\%$ of the acumulated $\sigma$. Therefore it is sufficient to take the reduced space of dimension $k = 1$. 
Note that in this 1D case, exploring the $\bz$ space consists in taking sets of $3$ neighbouring snapshots for the first level of connectivity,  $5$ snapshots for the second level, etc.

\begin{figure}[]
	\centering
	PCA \hspace{5.5cm} kPCA
	\includegraphics[width=0.9\textwidth]{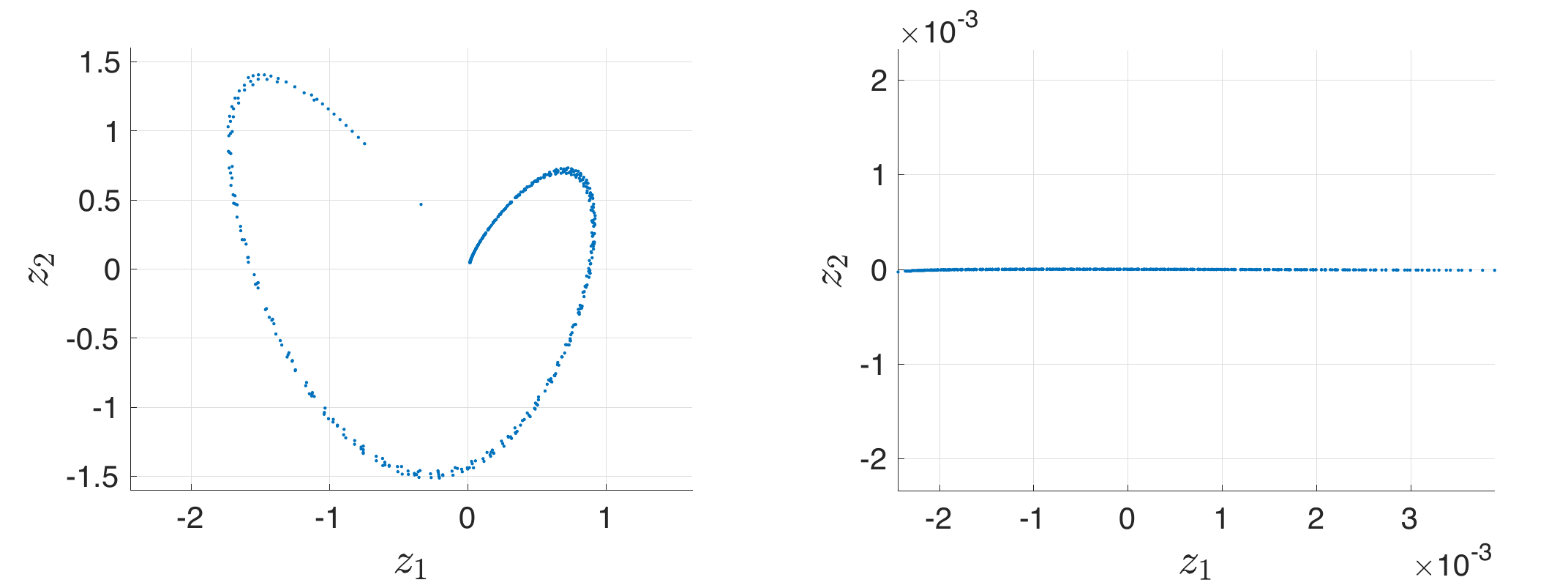}	 	
	\caption{\textit{1D problem.} Two first components of the snapshots in the reduced space from PCA and kPCA. PCA fails to identify the intrinsic dimension: $k=1$ contains $8.90\%$ of acumulated $\sigma$ and $k = 2$ contains $17.48\%$. On the other hand, kPCA captures $99.73\%$ of acumulated $\sigma$ with $k = 1$.}
	\label{fig:test1D-reducedspace}
\end{figure}

Figure \ref{fig:test1D-reducedspace} shows the $(z_1,z_2)$ components of the snapshots both in the PCA and the kPCA reduced spaces. It is clear how kPCA (with the proposed kernel) is able to identify the instrinsic dimension $k = 1$, being all the values of $z_{2}$ in a very narrow band. On the contrary, with PCA the range of the second component $z_{2}$ is similar to the range of $z_{1}$ (associated with only $8.90\%$ of the acumulated $\sigma$). 

\begin{figure}[]
	\centering
	\hspace{0.95cm}
	\hspace{-0.95cm}
	\includegraphics[width=0.485\textwidth]{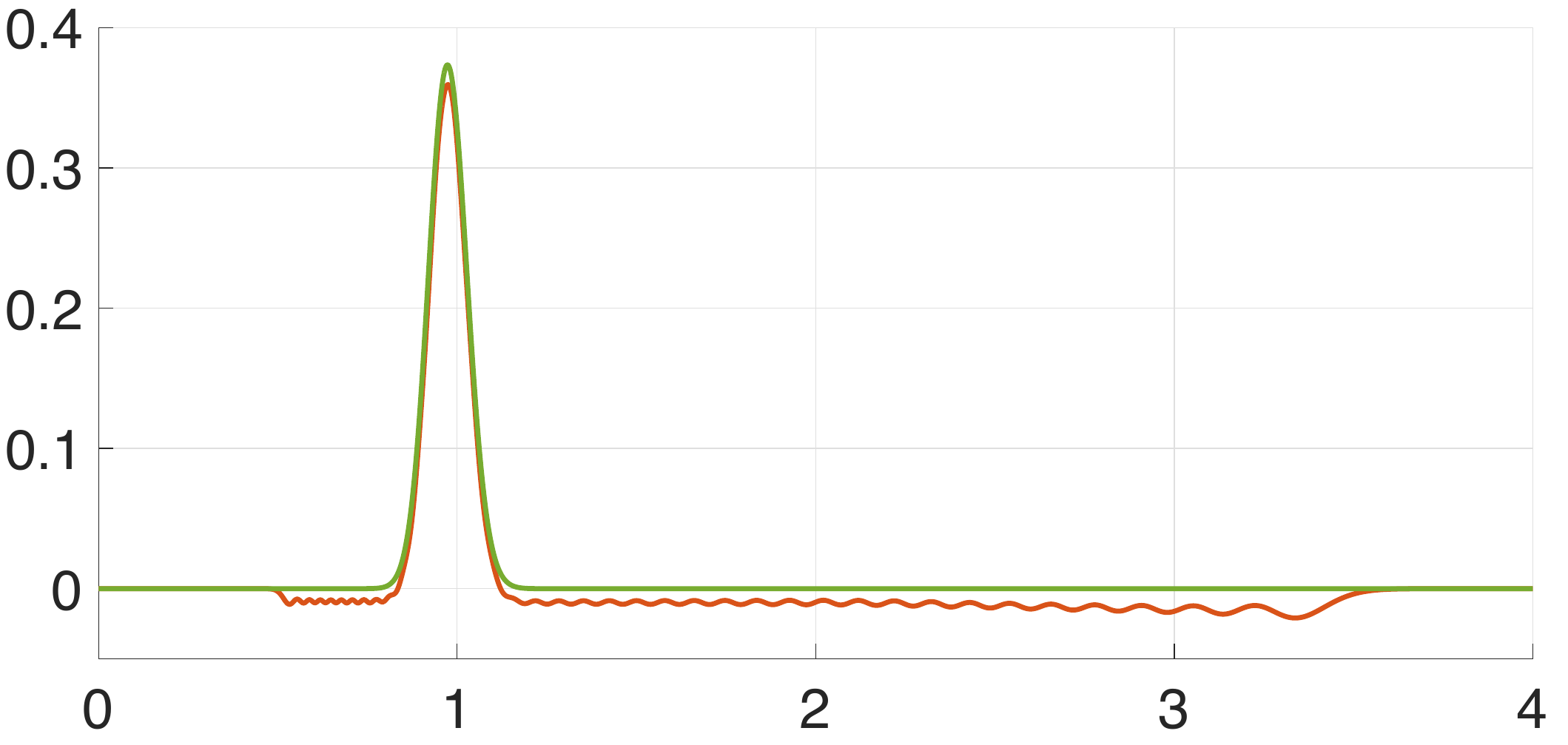}
	\hfill
	\includegraphics[width=0.485\textwidth]{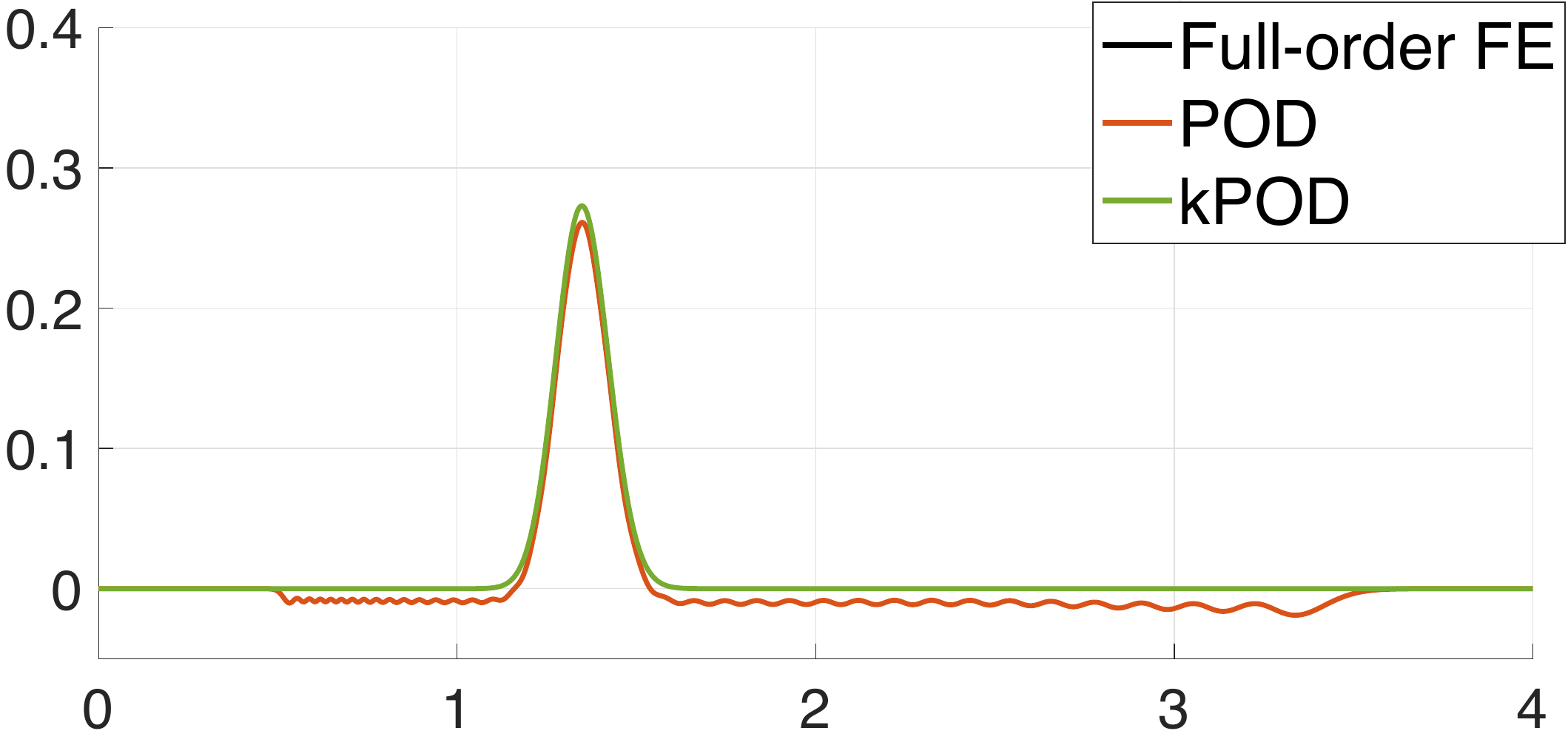}
	
	\raisebox{0cm}{\hspace{0.7cm} \text{$t = 0.25$} \hspace{5.7cm} \text{$t = 0.50$} }\hfill
	
	\vspace{3mm}
	\includegraphics[width=0.95cm]{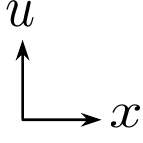}
	\hspace{-0.95cm}
	\raisebox{0.7cm}{\includegraphics[width=0.485\textwidth]{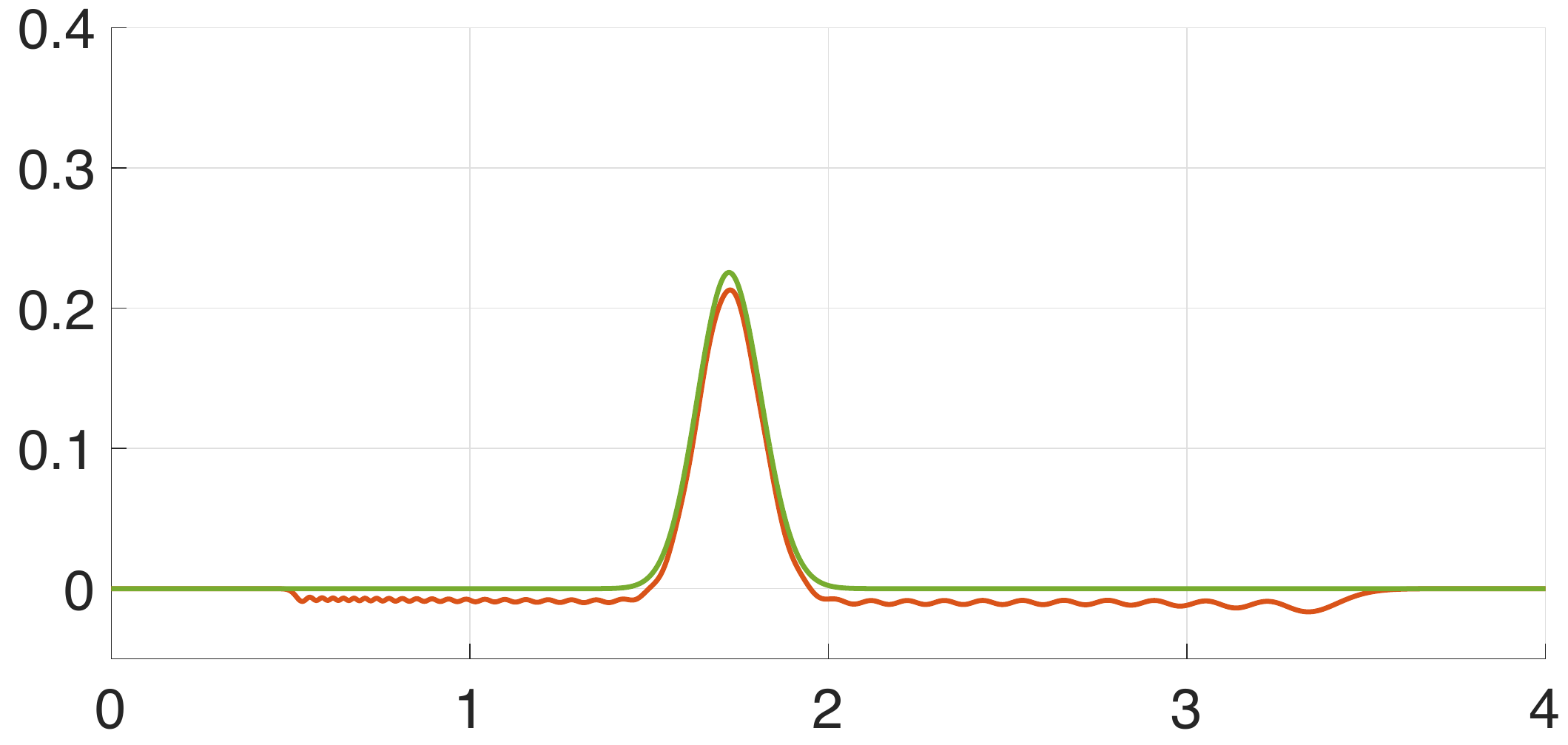}}
	\hfill
	\raisebox{0.7cm}{\includegraphics[width=0.485\textwidth]{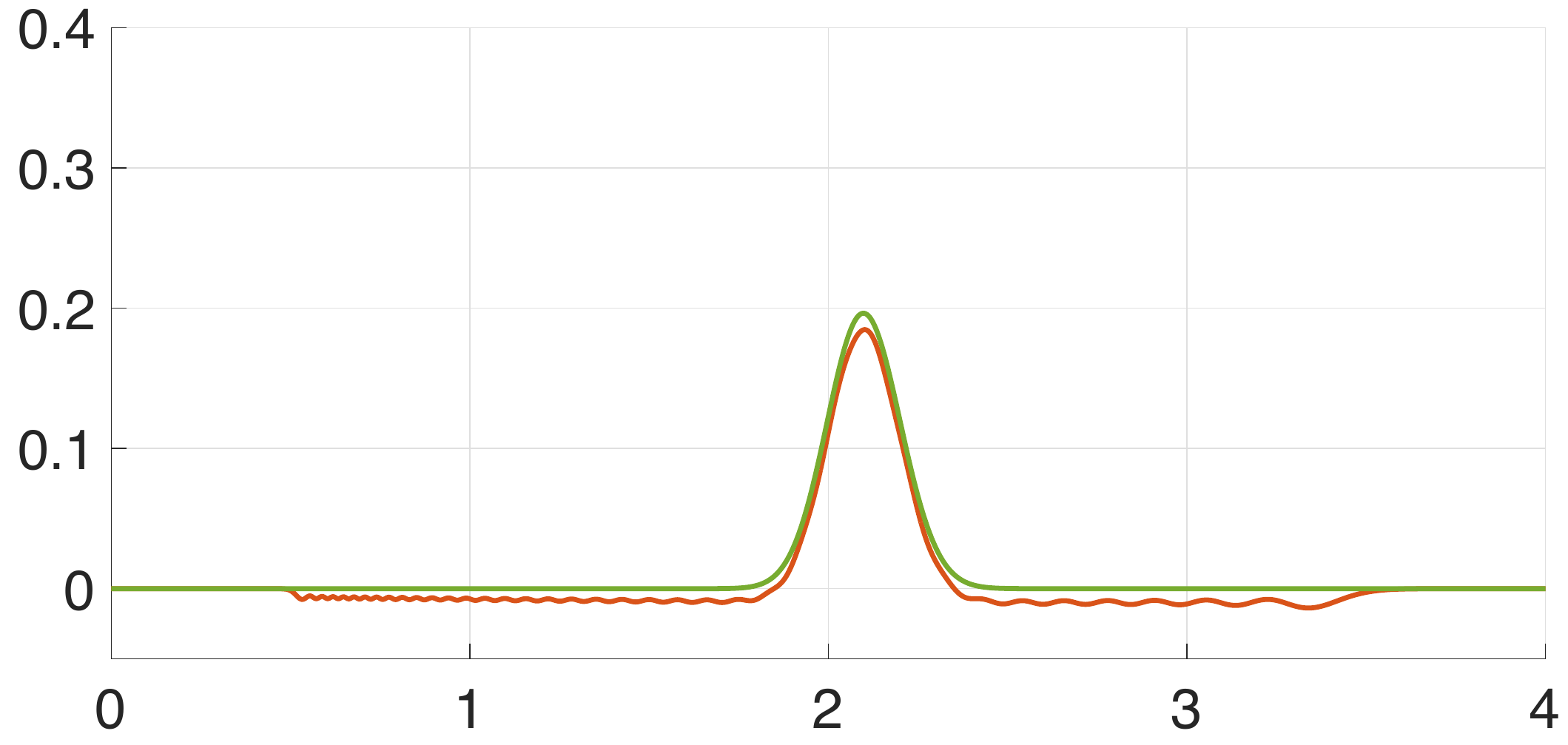}}
	
	\vspace{-0.7cm}
	\raisebox{0cm}{\hspace{0.7cm} \text{$t = 0.75$} \hspace{5.7cm} \text{$t = 1.00$} }\hfill
	
	\vspace{0.3cm}
	\caption{\textit{1D problem.} Reference full-order FE solution and corresponding approximations obtained with POD and kPOD (with $1/2$-connectivity levels) at different time steps. kPOD leads to more accurate solutions, undistinguishable from the reference solution in the plots.}
	\label{fig:test1D-plots}
\end{figure}

The accuracy of two formulations to obtain new solutions is analysed next.
We use the two methods to integrate in time for an advection velocity $v = 1.5$, which is not in the training set, taking increments in the temporal scheme of $\Delta t = 0.005$ until a final time $t = 1$.
The solutions at some intermediate time steps are shown in Figure \ref{fig:test1D-plots}. 

The POD solutions unphysically oscillate for this particular problem. Even though the obtained profile manages to follow the peak of the reference solution as time evolves, the relative errors with respect to the full-order solution are over $10^{-1}$, see Table \ref{table:test1D-POD}.
With kPOD, on the other hand, solutions do not oscillate and the relative errors decrease significantly. As listed in Table \ref{table:test1D-kPOD-steps}, errors for kPOD are of the order of $10^{-4}$.
Also, notice that the average number of steps in the optimal path is below $4$, with an average $\tilde{k}$ around $7$ during the whole process.

\begin{table}[]
	\centering
	\begin{tabular}{cc}
		\toprule
		Time  & Relative error \\ \midrule
		$0.25$  & $1.79\cdot 10^{-1}$ \\
		$0.50$ & $1.91\cdot 10^{-1}$ \\
		$0.75$ & $1.93\cdot 10^{-1}$ \\
		$1.00$ & $1.93\cdot 10^{-1}$ \\ \bottomrule
	\end{tabular}
	\caption{\textit{1D problem.} POD: relative errors at some time steps.} 
	\label{table:test1D-POD}
\end{table}

\begin{table}[]
	\centering
	\begin{tabular}{cccc}
		\toprule
		Time & Relative error & Steps opt path & Average $\tilde{k}$ \\ \midrule
		$0.25$  & $3.31\cdot 10^{-4}$  & $3.48$ & $7.00$ \\
		$0.50$ & $3.09\cdot 10^{-4}$  & $3.42$ & $6.99$ \\
		$0.75$ & $3.08\cdot 10^{-4}$  & $3.41$ & $7.05$ \\
		$1.00$ & $3.05\cdot 10^{-4}$  & $3.40$ & $7.04$  \\ \bottomrule
	\end{tabular}
	\caption{\textit{1D problem.} kPOD: relative errors, average number of steps in the optimal path and average $\tilde{k}$, for approximations with $1/2$-connectivity levels, at different time steps.}
	\label{table:test1D-kPOD-steps}
\end{table}

Higher accuracy is obtained with kPOD if increasing the levels of connectivity, accounting for more snapshots in the approximation of the solution. Table \ref{table:test1D-kPOD-levels} shows the errors at time $t=1$ for levels of connectivity $1$, $2$ and $3$ (corresponding to levels $2$, $3$ and $4$, respectively, in the final solution of Algorithm \ref{alg:OnlinePhase}). The results suggest that, as expected, $\tilde{k}$ increases with the level of connectivity, while the number of iterations in the optimal path tends to decrease.

\begin{table}[]
	\centering
	\begin{tabular}{cccc}
		\toprule
		Connectivity levels & Relative error & Steps opt path & Average $\tilde{k}$ \\ \midrule
		$1/2$ & $3.05\cdot 10^{-4}$ & $3.40$ & $7.04$ \\
		$2/3$ & $4.40\cdot 10^{-5}$ & $2.98$ & $14.11$\\
		$3/4$ & $3.95\cdot 10^{-7}$ & $2.98$ & $17.20$  \\ \bottomrule
	\end{tabular}
	\caption{\textit{1D problem.} kPOD: relative errors, average number of steps in the optimal path and average $\tilde{k}$ at time $t = 1.00$, for different levels of connectivity. }
	\label{table:test1D-kPOD-levels}
\end{table}

\subsection{Steady advection-diffusion problem in 2D}

The unknown function $u(x,y)$ takes values in the open domain $\Omega = ]-1,1[^2 \setminus B_{0.3}(0,0)$, where 
$B_{0.3}(0,0)=\left\{ (x,y)\in\RR^{2} \vert x^{2}+y^{2}\le 0.3^{2} \right\}$ is the circle or radius 0.3.
Let us consider the following steady advection-diffusion boundary value problem in $\Omega$
\begin{equation}\label{problem2D}
\left\{
\begin{aligned}
&-\bm{\nabla} \cdot \left( \nu \bm{\nabla} u \right) + \bv \bm{\nabla} u= 0  &&\text{in} &&\Omega, \\
&u = u_D  && \text{on} &&\Gamma_D = \partial\Omega \cap \{ x = -1 \} , \\
&\bm{\nabla} u \cdot \bn = 0  && \text{on} &&\partial\Omega \setminus \Gamma_D ,
\end{aligned}
\right.
\end{equation}
with diffusivity $\nu = 10^{-2}$ and $\bv(\bx)$ denoting the advection velocity field. The velocity field $\bv({\bx})$ has an average module of 10 and an inclination angle $\alpha\in[10^\circ,80^\circ]$. This is achieved by taking $\bv({\bx},\alpha) = 10\left(\bv_x({\bx}) \cos\alpha + \bv_y({\bx}) \sin\alpha\right)$, where velocities $\bv_x$ and $\bv_y$ correspond to unit horizontal and vertical flows in the domain. Velocity field $\bv_x$ (resp. $\bv_y$) is computed as deriving from a potential $\Phi$, and equal to $-\nabla \Phi$: potential $\Phi$ is the solution of a Laplace problem in $\Omega$, with boundary conditions $\nabla \Phi \cdot \bn = 1$ on the left (resp. top) side and $\nabla \Phi \cdot \bn = -1$ on the right (resp. bottom) side, complemented with homogeneous Neumann boundary conditions on the rest of $\partial \Omega$.

The Dirichlet condition in \eqref{problem2D} simulates a source (e.g. of concentration of some pollutant, or temperature) on the left portion of the boundary, with $u_D(x,y) = \frac{1}{\sigma\sqrt{2\pi}}\exp\left( \frac{-(y-\mu)^2}{2\sigma^2}\right)$, $\sigma = 2h$, $h = 0.02$ and with $\mu\in[-0.8,0.8]$ being the vertical coordinate of the source. 

The system is solved with the FEM on a triangular mesh of element size $h$, with
$20\,974$ elements and $10\,735$ nodes. The nodal FE solutions are vectors in $\RR^{10\,735}$.

\begin{figure}[]
	\centering
	\includegraphics[width=0.75\textwidth]{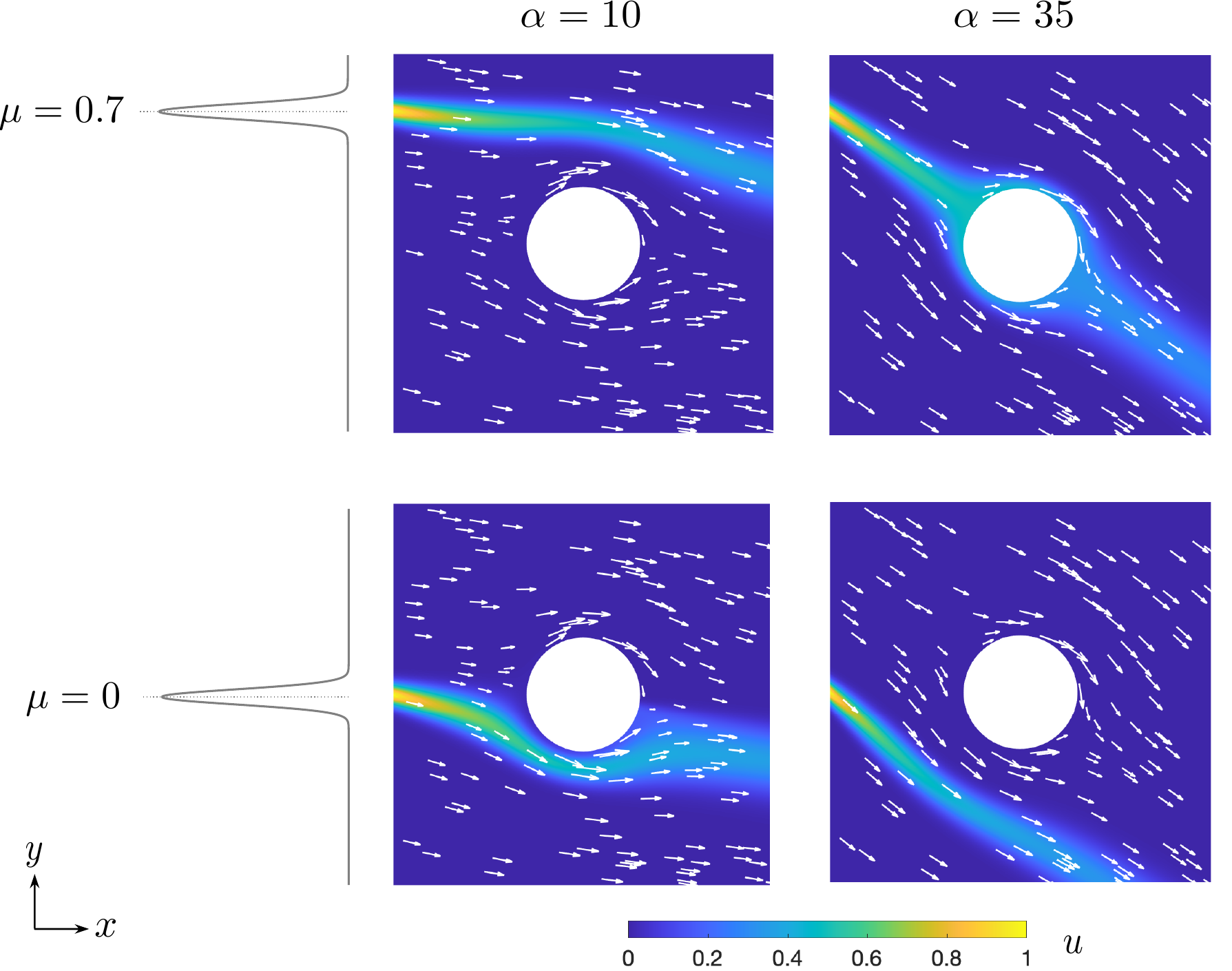}
	
	\caption{\textit{2D problem.} Some representative solutions for the system in \eqref{problem2D}. The pollution propagates above the island in $B_{0.3}(0,0)$ for $(\mu,\alpha) = (0.7,10)$, crosses the  island for $(0.7,35)$ and $(0,10)$, and propagates below the island for $(0,35)$. The white arrows indicate the direction of the advection velocity.}
	\label{fig:test2D-family}
\end{figure}

Figure \ref{fig:test2D-family} shows some representative solutions when varying the two parameters, $\alpha$ and $\mu$. The pollution propagates across the domain accordingly to $\mu$ and $\bv({\alpha})$, therefore it seems logical to think that solutions lie in a two-dimensional manifold. We can distinguish three different regimes depending on whether the wake of pollution is above, crossing or below the centered island in the domain. Crossing solutions are identified as those with $u > 0.1$ at some node on $\partial B_{0.3}(0,0)$. The lack of similarity between solutions in different regimes hints the need of a proper sampling  to obtain accurate results in all of them when using a reduced-order formulation.

\begin{figure}[]
	\centering
	\includegraphics[width = 0.9\textwidth]{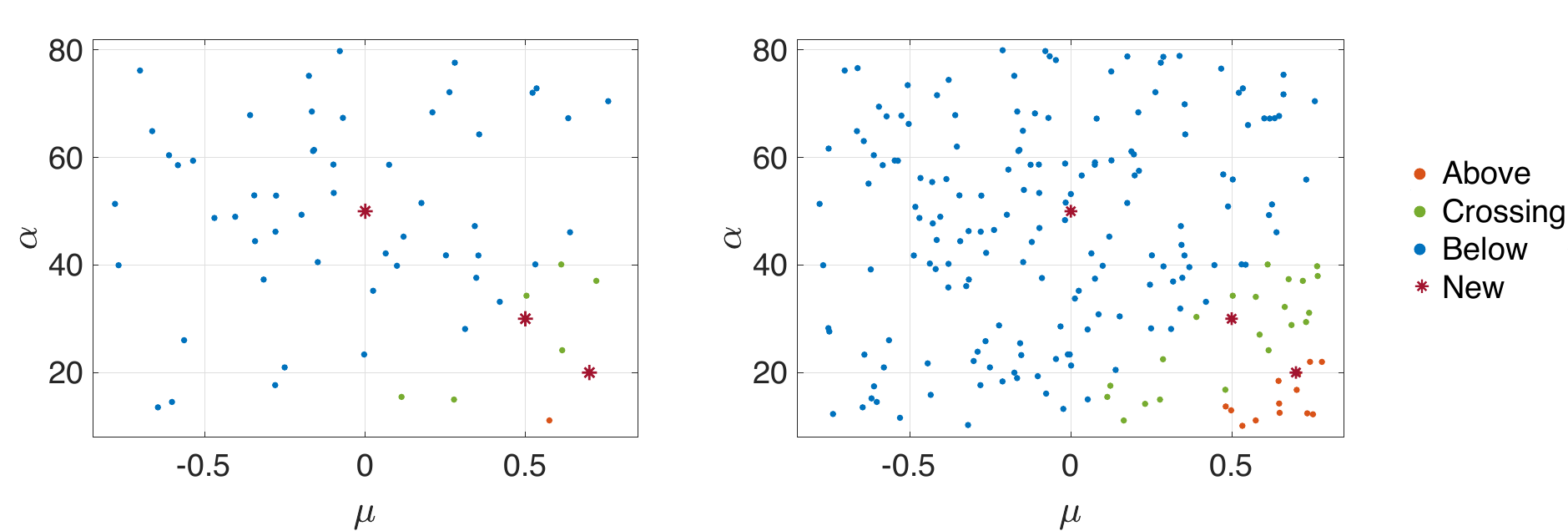}
	\hspace{2.7cm} \text{60 random samples} \hspace{2.3cm} \text{200 random samples} \hspace{1cm}
	
	\caption{\textit{2D problem.} Distribution of snapshots in the space of parameters $(\mu,\alpha)$. Different colors correspond to different regimes of the solutions (with paths below, crossing or above the island). Red asteriks indicate new configurations that are to be approximated with POD and kPOD.}
	\label{fig:test2D-params}
\end{figure}

We first consider a training set with $60$ snapshots corresponding to randomly distributed values for $(\mu,\alpha) \in [-0.8,0.8]\times[10,80]$, see Figure \ref{fig:test2D-params}. Note that most of the samples correspond to propagating paths below the island, while only one of them propagates above it. 

With a tolerance $\varepsilon = 10^{-8}$ in criterion \eqref{eq:CollectedVariance}, the PCA reduced space is of dimension $k = 58$. As in the previous example, PCA fails to identify the intrinsic dimension of the manifold, capturing $14.15\%$ of the acumulated $\sigma$ in the training set for $k = 2$.

\begin{figure}[]
	\centering
	\includegraphics[width=0.75\textwidth]{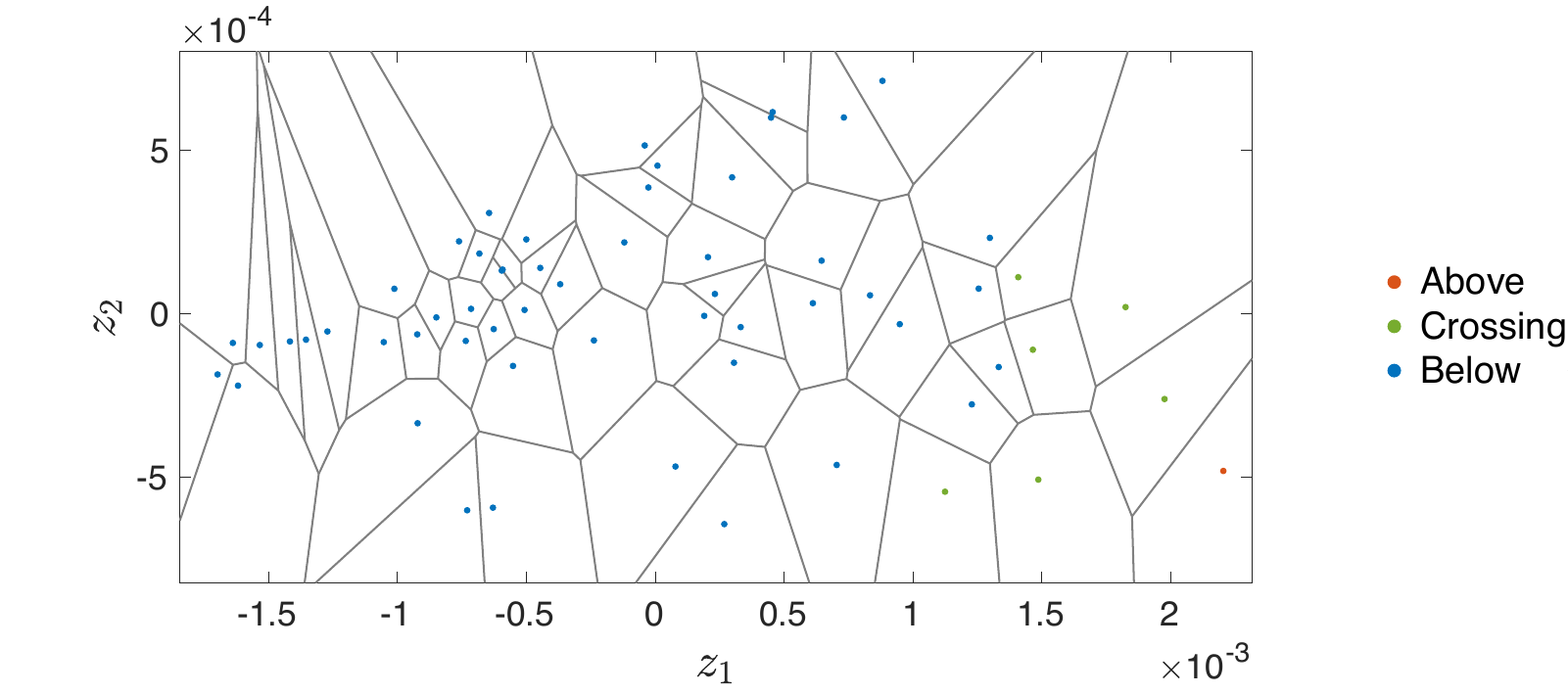}
	\caption{\textit{2D problem.} Voronoi tessellation in the kPCA reduced space of dimension $k = 2$, for the training set with $60$ random snapshots.}
	\label{fig:test2D-voronoi60}
\end{figure}

For kPOD, we define the kernel function
\begin{equation}\label{eq:kernel2D}
\kappa(\bu,\bv) = \exp \left(-\beta \left(\frac{ \| \bm{C}_{\Gamma_D}(\bu) - \bm{C}_{\Gamma_D}(\bv) \|^2}{2^2} +  \frac{\| \bm{C}_{out}(\bu) - \bm{C}_{out}(\bv)^j \|^2}{4^2} \right) \right),
\end{equation}
where $\bm{C}_{\Gamma_D}$ and $\bm{C}_{out}$ denote the centroid of the solution on $\Gamma_D$ and on the outlet $\partial\Omega\cap\left(\{y = -1\}\cup\{x= 1\} \right)$, respectively, as defined in \eqref{centroid_definition}. Taking $\beta = 10^{-3}$, the kernel explains $76.01\%$ of the acumulated $\sigma$ for $k = 1$, and $99.96\%$ for $k = 2$. We reduce to the intrinsic dimension $k = 2$, and the optimal path search is performed in the Voronoi tessellation in Figure \ref{fig:test2D-voronoi60}. As an intial approximation for the algorithm, we take $\tilde\bz_0 = F(\bu_{POD})$ (the image through the forward mapping of the POD solution).

In order to analyse the limit accuracy of kPOD with a quadratic approximation, we also consider the case in which all the snapshots in the domain are identified as neighbors. This is equivalent to POD with a quadratic approximation instead of the usual linear approach.
The solution from quadratic POD is expected to overcome the other ones in accuracy, since all the snapshots are involved. With the same tolerance $\varepsilon = 10^{-8}$, the resulting systems of equations have dimension $\tilde{k} = 1216$.

\begin{figure}[]
	\centering
	\hspace{2.5cm} POD \hspace{2.8cm} kPOD \hfill Quadratic POD \hspace{2cm}
	\vspace{1mm}
	
	\rotatebox{90}{{\hspace{3.2cm}}\rotatebox{-90}{c)} \hspace{3.4cm} \rotatebox{-90}{b)} \hspace{3.5cm} \rotatebox{-90}{a)}}
	\includegraphics[width = 0.88\textwidth]{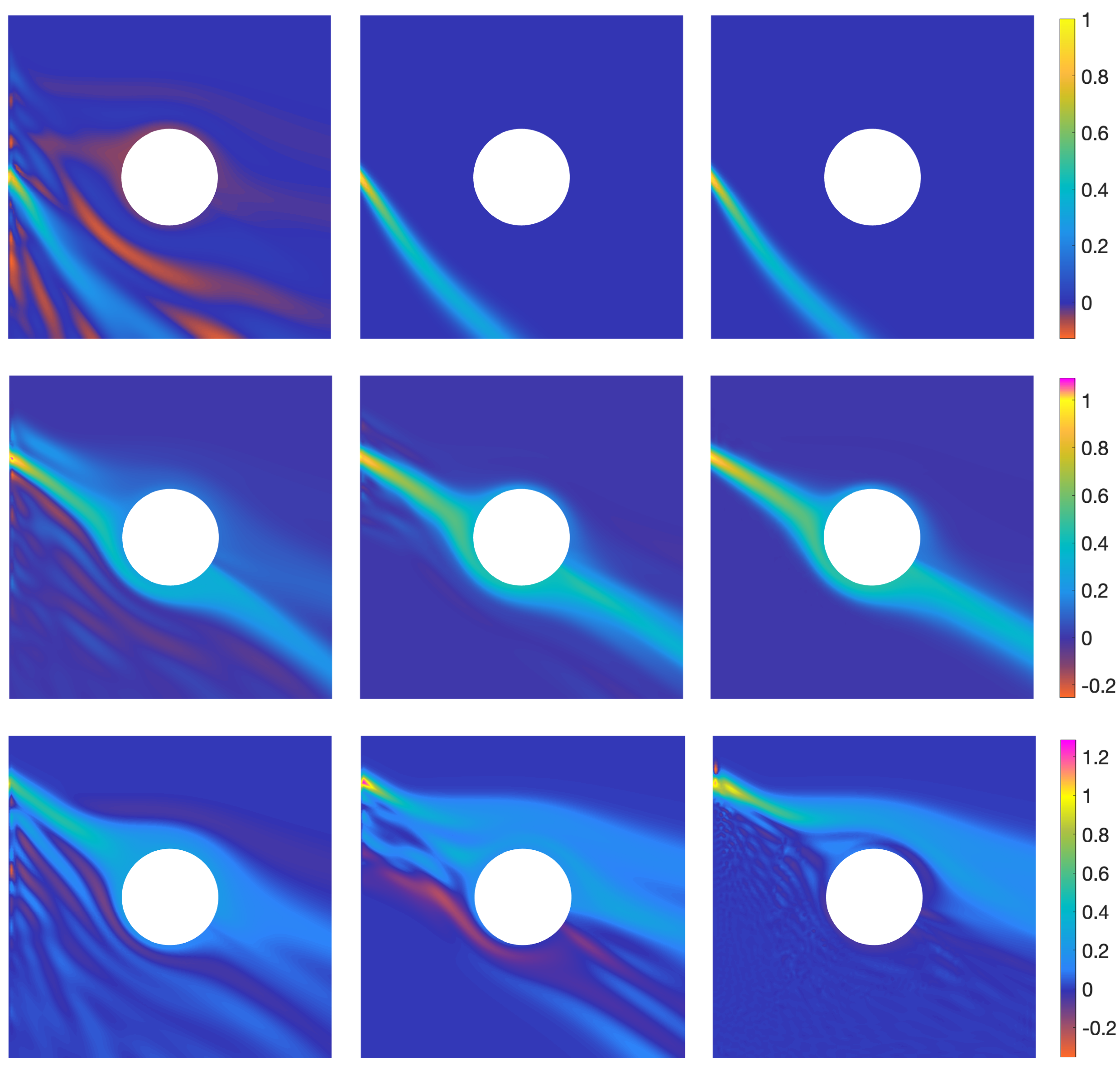}
	
	\caption{\textit{2D problem.} Solutions from POD, kPOD (with $1/2$-connectivity level) and quadratic POD, for a) $\mu = 0$, $\alpha = 50$ b) $\mu = 0.5$, $\alpha = 30$, and  c) $\mu = 0.7$, $\alpha = 20$, with a training set of $60$ random snapshots.}
	\label{fig:test2D-plots}
\end{figure}

Next, POD, kPOD and quadratic POD are used to approximate the solutions for three new points in the parameteric space, see Figure \ref{fig:test2D-params}. We take $(\mu,\alpha) = (0,50)$, $(0.5,30)$ and $(0.7,20)$, belonging to the below, crossing and above regimes, respectively. 
Figure \ref{fig:test2D-plots} shows the reduced-order solutions.

\begin{table}[]
	\centering
	\begin{tabular}{ll}
		\toprule
		Parameters & Relative error \\ \midrule
		$\mu = 0$, $\alpha = 50$ & $4.00\cdot 10^{-1}$ \\
		$\mu = 0.5$, $\alpha = 30$ & $4.21\cdot 10^{-1}$ \\ 
		$\mu = 0.7$, $\alpha = 20$ & $1.04$ \\\bottomrule
	\end{tabular}
	\caption{\textit{2D problem.} POD: relative errors for the new configurations of parameters, using the training set with $60$ snapshots.}
	\label{table:test2D-pod}
\end{table}

\begin{table}[]
	\centering
	\begin{tabular}{lllll}
		\toprule
		Parameters & Connectivity levels & Relative error & Steps opt path & $\tilde{k}$ \\ \midrule
		$\mu = 0$, $\alpha = 50$  & $1/2$ & $7.97\cdot 10^{-4}$ & $2$ & $43, 251$  \\
		& $2/3$ & $3.92\cdot 10^{-5}$ & $2$ & $251, 587$  \\ \cmidrule{2-5}
		$\mu = 0.5$, $\alpha = 30$  & $1/2$ & $9.29\cdot 10^{-2}$ & $2$ & $26,89$ \\
		& $2/3$ & $4.36\cdot 10^{-2}$ & $2$ & $89,169$ \\ \cmidrule{2-5}
		$\mu = 0.7$, $\alpha = 20$ & $1/2$ & $7.82\cdot 10^{-1}$ & $3$ & $13,13,43$ \\
		& $2/3$ & $4.97\cdot 10^{-1}$ & $4$ & $64,134,43,89$ \\ \bottomrule
	\end{tabular}
	\caption{\textit{2D problem.} kPOD: relative errors, number of steps in the optimal path and $\tilde{k}$ for each one of them, for the new configurations of parameters, using the training set with $60$ snapshots.}
	\label{table:test2D-kpod}
\end{table}

\begin{table}[]
	\centering
	\begin{tabular}{ll}
		\toprule
		Parameters & Relative error \\ \midrule
		$\mu = 0$, $\alpha = 50$ & $7.87\cdot 10^{-6}$ \\
		$\mu = 0.5$, $\alpha = 30$ & $1.91\cdot 10^{-2}$ \\ 
		$\mu = 0.7$, $\alpha = 20$ & $4.31\cdot 10^{-1}$ \\ \bottomrule
	\end{tabular}
	\caption{\textit{2D problem.} Quadratic POD: relative errors for the new configurations of parameters, using the training set with $60$ snapshots.}
	\label{table:test2D-podqua}
\end{table}

POD leads to results with oscillations and poor accuracy in the three cases. As listed in Table \ref{table:test2D-pod}, all
relative errors with respect to the FE reference solution are over $0.40$, with a relative error over $1$ for the solution corresponding to the above regime.

The accuracy in kPOD strongly depends on the number of snapshots in each one of the regimes, see Table \ref{table:test2D-kpod}. 
In the most-populated regime, the reduced solution for $(\mu,\alpha) = (0,50)$ has a relative error of $7.97\cdot 10^{-4}$ if taking $1/2$-connectivity levels to explore the reduced space, and the error decreases to $3.92\cdot 10^{-5}$ if $2/3$-connectivities are considered. 
Also, note that the optimal path only takes two iterations to converge.
For the solution in the crossing regime, the kPOD solution slightly oscillates, and oscillations are more evident for the new solution above the island, where the sampling is extremely poor, with only one snapshot. However, results are more accurate than those from POD in all cases. 

Results with quadratic POD, in Table \ref{table:test2D-podqua}, set the limit accuracy that we can reach with kPOD and a quadratic approximation for the given training set. Errors decrease and oscillations disappear (for the crossing-regime case) or become more subtle (for the above-regime one). Quadratic POD is an appealing option for its simple implementation and accuracy, but it is limited to datasets with a small number of snapshots. For a larger training set, for instance with $200$ snapshots, we may have memory limitations when performing the SVD of matrix $\bBL$ in \eqref{eq:localAver2}.

Finally, we repeat the POD and kPOD computations when enriching the training set up to $200$ snapshots, see Figure \ref{fig:test2D-params}. Again for $\varepsilon = 10^{-8}$, the reduced dimension of PCA is $k = 198$.
The kernel function \eqref{eq:kernel2D} enables kPCA to capture $99.97\%$ of the acumulated $\sigma$ for $k = 2$. 

\begin{figure}[]
	\centering
	\hspace{4cm} POD \hfill kPOD \hspace{5cm}
	\vspace{1mm}
	
	\rotatebox{90}{{\hspace{3.2cm}}\rotatebox{-90}{c)} \hspace{3.5cm} \rotatebox{-90}{b)} \hspace{3.6cm} \rotatebox{-90}{a)}}
	\includegraphics[width=0.65\textwidth]{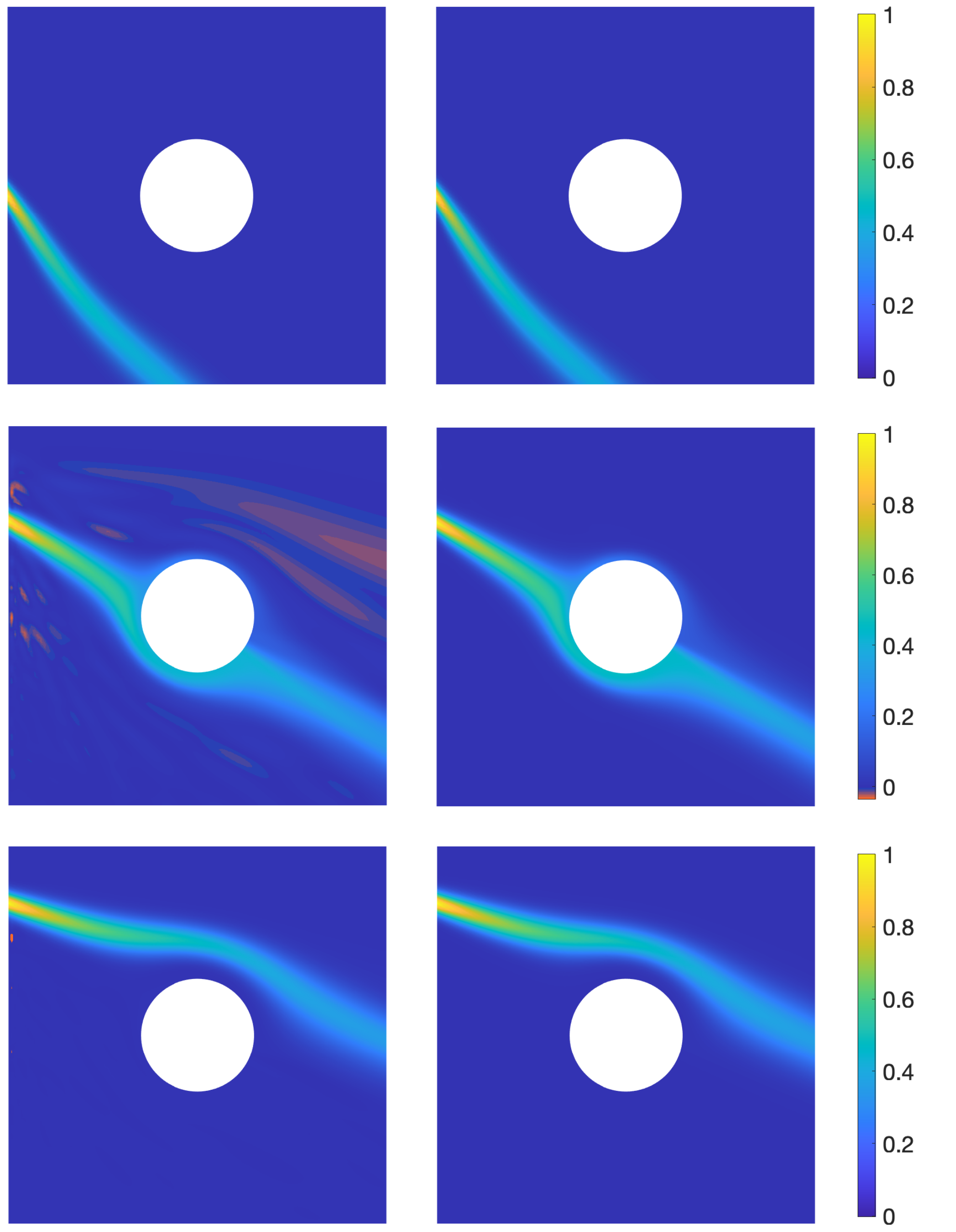}
	\caption{\textit{2D problem.} Solutions from POD and kPOD (with $1/2$-connectivity level), for a) $\mu = 0$, $\alpha = 50$ b) $\mu = 0.5$, $\alpha = 30$, and  c) $\mu = 0.7$, $\alpha = 20$, with a training set of $200$ random snapshots.}
	\label{fig:test2D-plots200}
\end{figure}

Figure \ref{fig:test2D-plots200} shows the POD and kPOD solutions for the same three cases as before. With POD, the errors decrease with the more populated training set, see Table \ref{table:test2D-pod200}, but we still observe oscillations in the tests on the crossing and above regimes.
kPOD leads to stable solutions in all cases. The errors decrease significantly, specially for the $2/3$-connectivities, see  Table \ref{table:test2D-kpod200}.

\begin{table}[h]
	\centering
	
	\begin{tabular}{ll}
		\toprule
		Parameters & Relative error \\ \midrule
		$\mu = 0$, $\alpha = 50$ & $4.83\cdot 10^{-3}$ \\
		$\mu = 0.5$, $\alpha = 30$ & $6.99\cdot 10^{-2}$ \\ 
		$\mu = 0.7$, $\alpha = 20$ & $1.42\cdot 10^{-2}$ \\ \bottomrule
	\end{tabular}
	\caption{\textit{2D problem.} POD: relative errors for the new configurations of parameters, using the training set with $200$ snapshots.}
	\label{table:test2D-pod200}
	
	\begin{tabular}{lllll}
		\toprule
		Parameters & Connectivity levels & Relative error & Steps opt path & $\tilde{k}$ \\ \midrule
		$\mu = 0$, $\alpha = 50$  & $1/2$ & $8.69\cdot 10^{-7}$ & $2$ & $76,433$  \\
		& $2/3$ & $1.58\cdot 10^{-8}$ & $2$ & $433,1232$  \\ \cmidrule{2-5}
		$\mu = 0.5$, $\alpha = 30$  & $1/2$ & $1.23\cdot 10^{-3}$ & $3$ & $26,26,251$ \\
		& $2/3$ & $6.45\cdot 10^{-5}$ & $3$ & $251,251,552$ \\ \cmidrule{2-5}
		$\mu = 0.7$, $\alpha = 20$ & $1/2$ & $5.70\cdot 10^{-4}$ & $2$ & $13,89$ \\
		& $2/3$ &  $9.09\cdot 10^{-5}$ & $2$ & $89,188$ \\ \bottomrule
	\end{tabular}
	\caption{\textit{2D problem.} kPOD: relative errors, number of steps in the optimal path and $\tilde{k}$ for each one of them, for the new configurations of parameters, using the training set with $200$ snapshots.}
	\label{table:test2D-kpod200}
\end{table}

\section{Concluding remarks}\label{sec:Conclu}
The numerical strategy presented in this paper, the kPOD, aims at alleviating the computational complexity of the parametric model by invoking nonlinear dimensionality-reduction techniques. More precisely, the kPCA is selected as an alternative to the linear PCA used in POD. 

The kPCA backward mapping (from the reduced space to the input space) is determined locally, using the snapshots surrounding the point to be mapped. This is inducing a local character of the kPOD algorithm, and has suggested devising an iterative exploration of the reduced space, based on a Delaunay tessellation of the cloud of snapshots (in the low-dimensional space).

Another novel idea introduced here is to increase the basis of snapshots in the input space by including the cross-products of the snapshots, that is, the quadratic terms. This is introducing new features of the solutions that contribute to a better resolution in the kPOD. In practice, this is only affordable in a kPOD setting where the solutions are computed locally (with a reduced number of neighbouring snapshots). The number of resulting quadratic terms arising from a global representative family of snapshots is too large and compromises the use of the standard POD.

Moreover, the selection of the kernel to define the kPCA reduction is carried out using considerations based on the nature of the problem. The kernel is taken such that it extracts features of the solution that discriminate the difference in the associated parameters (e.g. the location of the barycenter of the distribution as an indicator of the time and the diffusion).

The combination of the different novel elements presented in the paper is increasing the portfolio of numerical tools to devise reduced-order strategies dealing with complex problems. Taking advantage of the nonlinear manifold allows increasing the accuracy of the solution with respect to the standard POD, noting that this methodology does not necessarily reduce the computational cost.

The performance of the presented strategies in two parametric advection-diffusion problems (one transient, one steady-state) shows a considerable improvement with respect to the standard POD.

\section*{Acknowledgements}
This work is partially funded by Spanish Ministry of Economy and Competitiveness and Generalitat de Catalunya (DPI2017-85139-C2-2-R and 2017-SGR-1278)

\bibliographystyle{plain} 
\bibliography{Ref}

\end{document}